\documentclass[11pt]{article}
\usepackage{amssymb}
\usepackage{url}
\usepackage[dvips]{graphicx}
\usepackage{graphics}
\usepackage[dvips]{color}
\usepackage{xcolor}
\usepackage[figuresright]{rotating}
\usepackage{multirow}
\usepackage[top=2cm, bottom=2cm, left=2cm, right=2cm]{geometry}
\usepackage{algorithm}
\usepackage{algorithmicx}
\usepackage{algpseudocode}
\usepackage{algpseudocode}
\usepackage{amsmath}
\usepackage{amsthm}
\usepackage{amsmath}
\usepackage{amsfonts}
\usepackage{epic,eepic}
\usepackage{epsfig}
\usepackage{booktabs} 
\usepackage[toc,page,title,titletoc,header]{appendix}
\usepackage[font={small}]{caption}
\newtheorem{thm}{Theorem}[section]
\newtheorem{lem}[thm]{Lemma}
\newtheorem{cor}[thm]{Corollary}

\newtheorem{rem}[thm]{Remark}

\def\dg{{\rm diag}}
\def\Dg{{\rm Diag}}

\def\ra{{\rm Range}}
\def\rk{{\rm rank}}
\def\proof {\pn {Proof.} }

\def\endproof{\hfill $\Box$ \vskip .5cm}

\newcommand\norm[1]{\left\lVert#1\right\rVert}

\usepackage{amsfonts}
\def\Dg{{\rm Diag}}
\date{\today}
\author{Rujun Jiang\thanks{School of Data Science, Fudan University, China, rjjiang@fudan.edu.cn}
\and Duan Li\thanks{Corresponding author. Department of Systems Engineering and Engineering Management, The Chinese University of Hong Kong, Hong Kong,  dli@se.cuhk.edu.hk}}

\title{Second Order Cone Constrained Convex Relaxations for
Nonconvex Quadratically Constrained Quadratic
Programming}
\usepackage{varioref}
\def\endproof{\hfill $\Box$ \vskip .7cm}
\def\proof {\textit{{Proof.}} }

\begin{document}
\maketitle
\begin{abstract}

In this paper, we present new convex relaxations for  nonconvex quadratically constrained quadratic programming (QCQP) problems. While recent research has focused on strengthening
convex relaxations using reformulation-linearization technique (RLT), the state-of-the-art methods lose their effectiveness when dealing with (multiple) nonconvex quadratic constraints in QCQP.
In this research, we decompose and relax each nonconvex constraint to two second order cone (SOC) constraints and then linearize the products
of the SOC constraints and linear constraints to construct some effective new valid constraints.
Moreover, we extend the reach of the RLT-like techniques for almost all different types of constraint-pairs (including valid inequalities by linearizing the product of a pair of SOC constraints, and the
Hadamard product or the
Kronecker product of two respective valid linear matrix inequalities), examine dominance relationships among different valid inequalities, and explore almost all possibilities of gaining benefits from generating valid constraints. Especially, we successfully demonstrate that applying RLT-like techniques to additional redundant linear constraints could reduce the relaxation gap significantly. We demonstrate the efficiency of our results
with numerical experiments.
\end{abstract}

\section{Introduction.}
We consider in this paper the following class of quadratically constrained quadratic programming (QCQP) problems:
\begin{eqnarray*}
 \rm{(P)}~~~
&~  \min& x^TQ_0x+c_0^Tx\notag \\
& ~ \rm{s.t.}&x^TQ_ix+c_i^Tx+d_i\leq0,~i=1,\ldots,l,\\
&~&a_j^Tx\leq b_j,~j=1,\ldots,m,
\end{eqnarray*}
where $Q_i$ is an $n\times n$ symmetric matrix, $c_{i}\in \Re^n$, $i=0\ldots,l$, $d_i\in \Re$, $i=1\ldots,l$ and  $a_j\in \Re^n,~b_j\in \Re$, $j=1, \ldots, m$.
Without loss of generality, we assume that
$Q_i$ is not a zero matrix for $i=1,\ldots,l$. We further partition the quadratic constraints into the following two groups:
\begin{eqnarray*}&\mathcal{C}&=\{i:~Q_i \text{ is positive semidefinite},~i=1,\ldots,l\}, \\
&\mathcal{N}&=\{i:~Q_i \text{ is not positive semidefinite},~i=1,\ldots,l\}.
\end{eqnarray*}
Without loss of generality, we assume in this paper the cardinality of $\mathcal{C}$ is $k$ ($k\leq l$). QCQP problems arise in various areas, for example, combinatorial optimization, portfolio selection problems, economic equilibria, 0--1 integer programming and various applications in engineering. In the past few decades, QCQP has been  widely investigated in the
literature (see, e.g., \cite{anstreicher2012convex,ben2014hidden,xia2016s,bao2011semidefinite,fujie1997semidefinite,kim2016lagrangian,linderoth2005simplicial,ye2003new,zheng2011convex}), due to its elegance in formulation and a wide spectra of applications.

QCQP in general is NP-har \cite{sahni1974computationally}, even when it only has  linear constraints \cite{pardalos1991quadratic}, although some special cases of QCQP are polynomially solvable
\cite{anstreicher1999strong,anstreicher2000lagrangian,burer2013second,sturm2003cones,beck2006strong}. As a global optimal solution of QCQP is generally hard to compute due to its NP-hardness, based on various kinds of relaxations,
branch and bound methods  have been developed in the literature to find exact solutions for QCQP problems; see, e.g.,
\cite{linderoth2005simplicial,burer2008finite}.
It is well known that the  efficiency of a branch and bound method depends on two major factors: the quality of the relaxation bound and its associated
computational cost.
Recent decades have witnessed an increasing attention on constructing convex relaxations enhanced with various valid inequalities.
The survey paper \cite{bao2011semidefinite} compared the computational speed and quality of the gaps of various semidefinite  programming (SDP) relaxations with different valid inequalities for QCQP problems.
Sherali and Adams \cite{sherali2013reformulation} first introduced the concept of ``reformulation-linearization technique"  (RLT) to
achieve a lower bound of problem (P). Anstreicher in \cite{anstreicher2009semidefinite}  proposed a theoretical analysis for successfully applying RLT constraints to remove a large portion of the feasible region for the relaxation, and
suggested that a combination of SDP and RLT constraints leads to a tighter bound. This standpoint holds true for the  relaxations
with all other valid inequalities based on the idea behind RLT in this paper.
Sturm and Zhang \cite{sturm2003cones} developed the so-called SOC-RLT constraints (or called \textit{rank-2 second-order inequalities} in  \cite{yang2013two,zheng2011convex}) to solve the problem of minimizing a quadratic objective function subject to a convex quadratic constraint and a linear constraint exactly when combined with its SDP relaxation. More specifically, they rewrote a convex quadratic constraint as a second order cone (SOC) constraint and linearized the product of the SOC  and linear constraints.
Burer and Saxena \cite{burer2012milp} discussed how to utilize the SOC-RLT constraints to get a tighter bound than the SDP+RLT relaxation for general mixed integer QCQP problems.
Recently, Burer and Yang \cite{burer2013trust} demonstrated that the SDP+RLT+(SOC-RLT) relaxation has no gap in an extended trust region problem of minimizing a quadratic function subject to a unit ball and multiple linear constraints, where the linear constraints do not intersect with each other in the interior of the ball.

However, all methods mentioned above lose their effectiveness when dealing with (multiple) nonconvex quadratic constraints in QCQP problems. The  state-of-the-art in dealing with nonconvex quadratic constraints is to directly lift the quadratic terms as the basic SDP relaxation does.
This recognition and the success of combining SDP relaxations with RLT and SOC-RLT constraints (for convex quadratic constraints) motivate our study in this paper. Using the basic ideas behind SOC-RLT constraints, our method constructs valid inequalities based on linearizing the product of  the nonconvex
quadratic constraints and linear constraints, and  performs better than the state-of-the-art convex relaxations for problem (P).
We call our newly developed valid  inequalities  Generalized SOC-RLT (GSRT) constraints.
For simplicity of analysis, we call any nonconvex quadratic constraint type-A and a nonconvex quadratic constraint  $x^TQ_ix+c_i^Tx+d_i\leq0$ type-B if $c_i\in \ra(Q_i)$.
 Although this range condition for type-B could be  numerically hard to check in general, it  can be readily verified in some special cases, e.g., $Q_i$ is nonsingular or  $c_i\in \ra(Q_i)$ is known to hold in advance for some specific problem data set. To construct  GSRT constraints, we first introduce a new augmented variable  $z_i$
corresponding to each nonconvex constraint $x^TQ_ix+c_i^Tx+d_i\leq0,~i\in\mathcal{N}, $   and then decompose the matrix $Q_i$ according to the signs of its eigenvalues such that $Q_i=L_i^TL_i-M_i^TM_i$.
Depending on different techniques in handling the linear term, the decomposition of $x^TQ_ix+c_i^Tx+d_i\leq0$ further results in
two types of GSRT, i.e., type-A GSRT constraint (GSRT-A) and type-B GSRT constraint (GSRT-B) as follows.
 GSRT-A is derived from the equivalence between $x^TQ_ix+c_i^Tx+d_i\leq0$ and the following two constraints,
 \begin{eqnarray}
&&\norm{\left.\left(\begin{array}{c}L_ix\\\frac{1}{2}(c_i^Tx+d_i+1)\end{array}\right)\right.}\leq z_{i},\label{GSRT-A1}\\
&&\norm{\left(\begin{array}{c}M_ix\\\frac{1}{2}(c_i^Tx+d_i-1)\end{array}\right)}= z_{i},\label{GSRT-A2}
 \end{eqnarray}
where the equivalence is easily derived by substituting (\ref{GSRT-A2}) into (\ref{GSRT-A1}). If a type-B quadratic constraint holds for index $i$ with $c_i\in{\rm Range}(Q_i)$, GSRT-B constraints are then constructed by decomposing $x^TQ_ix+c_i^Tx+d_i\leq0$ in one of the following two different  ways:
\begin{itemize}
\item i)
if $\frac{1}{4}(c_i^TQ_i^{\dagger}c_i)-d_i\geq0$, we decompose $x^TQ_ix+c_i^Tx+d_i\leq0 $ as
\begin{eqnarray}
\left.\begin{array}{l}\norm{L_i(x+x_{0})}\leq z_{i},\\
\norm{\left(\begin{array}{c}M_i(x+x_{0})\\\Delta\end{array}\right)}= z_{i},
\end{array}\right.\label{b1}
\end{eqnarray}
where $\Delta=\sqrt{\frac{1}{4}(c_i^TQ_i^{\dagger}c_i)-d_i}$,
  $x_0=\frac{1}{2}Q_i^{\dagger}c_i$ and $A^{\dagger}$ denotes the \textit{Moore--Penrose} pseudoinverse for matrix $A$.
\item ii) if $\frac{1}{4}(c_i^TQ_i^{\dagger}c_i)-d_i<0$, we decompose $x^TQ_ix+c_i^Tx+d_i\leq0 $ as
\begin{eqnarray}
\left.\begin{array}{l}\norm{\left(\begin{array}{c}L_i(x+x_{0})\\\Delta\end{array}\right)}\leq z_{i}, \\
\norm{M_i(x+x_{0})}= z_{i},
\end{array}\right.\label{b2}
\end{eqnarray}
 where  $\Delta=\sqrt{d_i-\frac{1}{4}(c_i^TQ_i^{\dagger}c_i)}$,  $x_0=\frac{1}{2}Q_i^{\dagger}c_i$.
\end{itemize}
 Since  the   equality constraint  (\ref{GSRT-A2}) is nonconvex and   intractable, we relax\ (\ref{GSRT-A2}) to an inequality  to obtain an SOC constraint (which is convex and tractable),
\begin{equation}
 \norm{\left(\begin{array}{c}M_ix\\\frac{1}{2}(c_i^Tx+d_i-1)\end{array}\right)}\leq z_{i}.\label{GSRT-A2i}
\end{equation}
Multiplying any linear constraint to both sides of the above two kinds of SOC constraints in (\ref{GSRT-A1}) and (\ref{GSRT-A2i}), respectively, and linearizing the products lead to additional valid inequalities.
Moreover, we construct   valid equalities by linearizing the squared form of  (\ref{GSRT-A2}), i.e., linearizing the following equality,
\begin{equation}
xM_i^TM_ix+\frac{1}{4}(c_i^Tx+d_i-1)^{2}=z_{i}^2.\label{squred}
\end{equation}
The GSRT-A constraints consist of SOC constraints in  (\ref{GSRT-A1}) and (\ref{GSRT-A2i}),  the linearization of the products of SOC constraints in (\ref{GSRT-A1}) and (\ref{GSRT-A2i}) with
any original linear constraint, and  the linearization of (\ref{squred}).
With similar techniques, we can construct GSRT-B constraints according to the different decomposition schemes of $x^TQ_ix+c_i^Tx+d_i\leq0$, given in (\ref{b1}) and (\ref{b2}), respectively.
Note that GSRT-A constraints can be generated from any pair of a nonconvex quadratic constraint and a linear constraint, but GSRT-B constraints can only be generated from those pairs under the range condition $c_i\in{\rm Range}(Q_i)$. That is, we can always construct GSRT-A, but have limited ability to construct GSRT-B only under the range condition $c_i\in{\rm Range}(Q_i)$.
We then prove that the GSRT relaxation, which stands for the SDP relaxation enhanced  with RLT, SOC-RLT and GSRT constraints, achieves a much tighter
lower bound for problem (P) than the sate-of-the-art relaxation in the literature.

Another RLT-based technique in the literature is to introduce and attach additional redundant linear constraints to the original QCQP problem and then apply  the RLT and SOC-RLT techniques.
Zheng et al.  \cite{zheng2011convex} proposed  a decomposition-approximation method for generating convex relaxations to get a tighter lower bound than the SDP+RLT+(SOC-RLT) bound.
Enlightened by  the decomposition-approximation method in \cite{zheng2011convex}, we introduce a new relaxation by generating
extra RLT, SOC-RLT and GSRT constraints
with extra redundant linear inequalities.  We further demonstrate that this relaxation dominates  the decomposition-approximation method in \cite{zheng2011convex}
  for  problem (P)  with an extra nonnegativity constraint $x\geq0$.

Inspired by the GSRT constraints, we also explore and construct a new class of valid inequalities by linearizing the product of any pair of SOC constraints,
termed SOC-SOC-RLT (SST) constraint.
Moreover, we demonstrate that
 this new class of valid inequalities is equivalent to a valid  linear matrix inequality (LMI) formed by a submatrix of the Kronecker product constraint proposed
 in  \cite{anstreicher2017kronecker}, termed Kronecker SOC-RLT (KSOC) constraint. However, as the KSOC constraint is a large-scale  LMI, its dimensionality may prevent its direct application from practical implementation. We thus discuss the tradeoff between using KSOC and its submatrices with respect to the bound quality and computational costs. We also investigate several other KSOC constraints and their dominance relationship with the valid inequalities discussed in this paper.

 We illustrate below the different kinds of valid inequalities generated by RLT-like technique, i.e., linearizing the product of the left hand side yields the
valid inequalities on the right hand side, and also indicate in the list the sections (or subsections) in which different RLT-like techniques are developed,
$$\begin{array}{rcll}\rm L\times L&\Longrightarrow&\rm  RLT ~([22])&\text{Section 2.1,}\\
\rm SOC(convex)\times L&\Longrightarrow &\rm SOC\mbox{-}RLT~([25])&\text{Section 2.1,}\\
\rm SOC(nonconvex)\times L&\Longrightarrow &\rm GSRT&\text{Section 2.2,}\\
\rm M(\succeq0)\circ M(\succeq0)&\Longrightarrow &\rm HSOC~([29])~~~~~~ &\text{Section 3,}\\
\rm SOC\times SOC&\Longrightarrow &\rm SST~~~~~~~~~~~~~~~~~&\text{Section 4,}\\
\rm M(\succeq0)\otimes M(\succeq0)&\Longrightarrow &\rm KSOC(\succeq0)~([3] \text{ and this paper})~&\text{Section 5,}
\end{array}$$
where $\rm L$ represents a linear inequality constraint, $\rm SOC(convex)$ ($\rm SOC(nonconvex)$, respectively) represents an SOC constraint generated from a convex (nonconvex, respectively) constraint, $\rm M(\succeq0)$ represents an LMI, HSOC represents the valid inequalities generated by linearizing the Hadamard product of two valid LMIs (expressed in (\ref{hada}) later in the paper) in \cite{zheng2011convex} and KSOC represents the valid inequalities generated by linearizing the Kronecker product of two valid LMIs first derived in \cite{anstreicher2017kronecker}.

In general, there is no dominance relationship among
the valid inequalities RLT, SOC-RLT, GSRT and KSOC. Furthermore, although SST, HSOC and the valid LMI  given in (\ref{subksrt}) later in the paper are not dominated by RLT, SOC-RLT and GSRT, they are all dominated by a KSOC valid inequality as we  will prove in Section 5. When a new valid inequality has no dominance relationship with the existing constraints in the formulation, adding this additional valid inequality to the constraints should yield a tighter relaxation. So the guiding principle of our research is to extend the RLT-like technique to derive effective
valid inequalities to strengthen the SDP relaxation, especially to develop effective
valid inequalities from nonconvex quadratic constraints.

We summarize now the main contributions of this paper in the following three aspects.
\begin{itemize}
\item We derive the GSRT constraints, which represent the first attempt in the literature to construct new valid inequalities  for nonconvex quadratic constraints using RLT-like techniques.
\item We extend the reach of the RLT-like techniques for almost all different types of constraint-pairs and explore almost all possibilities of gaining benefits from generating valid constraints. We also successfully demonstrate that applying RLT-like techniques to additional redundant linear constraints could reduce the relaxation gap.
\item We examine possible dominance relationships among different valid inequalities generated from various RLT-like techniques. We also discuss the tradeoff between the tightness of the bound and the computational cost.
\end{itemize}

 The rest of the paper is organized as follows. In Section 2, we first review existing convex relaxations with various valid inequalities in
the literature and then propose our novel GSRT constraints. In Section 3, we apply RLT-like techniques to additional redundant linear constraint and demonstrate a dominance relationship of our method over the method in \cite{zheng2011convex}. We propose in Section 4 another class of valid inequalities, SST constraints, by
linearizing the product of two SOC constraints. In Section 5, we introduce KSOC constraints in the recent literature and show their relationships with the previous constraints discussed in the paper.
After we demonstrate good performance of GSRT from numerical tests in Section 6, we
offer our concluding remarks in Section 7.

\textbf{Notation}
We use $v(\cdot)$ to denote the optimal value of problem $(\cdot)$. Let $\norm{x}$  denote the Euclidean norm of $x$, i.e., $\norm{x}=\sqrt{x^Tx}$,
and $\norm{A}_F$  denote the Frobenius norm of a matrix $A$, i.e., $\norm{A}_F=\sqrt{tr(A^TA)}$.
The notation  $A\succeq0$ refers that  matrix  $A$ is positive semidefinite and the notation $A\succeq B$ implies that $A-B\succeq0$.
 The inner product of two symmetric
matrices is defined by $A\cdot B=\sum_{i,j=1,\ldots,n} A_{ij}B_{ij}$, where $A_{ij}$ and $B_{ij}$ are the $(i,j)$ entries of $A$ and $B$, respectively. We also use $A_{i,\cdot}$ and $A_{\cdot,i}$ to denote
the $i$th row and column of  matrix $A$, respectively. Notation $\rk(A)$ denotes the rank of matrix $A$.
We use $\dg(v)$, where $v$ is a column vector, to denote a diagonal matrix with its $i$th diagonal entry being $v_{i}$ and $\Dg(A)$ to denote
the column vector with its $i$th entry being $A_{ii}$. For a positive semidefinite matrix $A$  with spectral decomposition $A=U^TDU$, where $D$ is a diagonal matrix,
we use notation $A^{\frac{1}{2}}$ to denote $U^TD^{\frac{1}{2}}U$, where $D^{\frac{1}{2}}$ is a diagonal matrix with $\sqrt{D_{ii}}$ being its $i$th entry.

\section{Generalized SOC-RLT constraints}
In this section, we first present the basic SDP relaxation for problem (P) and its strengthened variants with RLT and SOC-RLT constraints
in the literature and then propose the new GSRT constraints.
\subsection{Preliminary}
Let us now first review some existing relaxations for problem $\rm(P)$ in the literature.
By lifting $x$ to  matrix  $X=xx^T$ and relaxing $X=xx^T$  to $X\succeq xx^T$,  which is further equivalent to $\left(
\begin{array}{cc}
1&x^T\\
x&X
\end{array}
\right)\succeq0$  due to the Schur complement, we have the following basic SDP relaxation for problem (P):
\begin{eqnarray}
 {\rm (SDP)}
&~\min &Q_0\cdot X+c_0^Tx\notag \\
&~\text{s.t.}& Q_i\cdot X+c_i^Tx+d_i\leq0,~i=1,\ldots,l,\label{QC}\\
&~&a_j^Tx\leq b_j,~j=1,\ldots,m,\label{LC}\\
&~&\left(
\begin{array}{cc}
1&x^T\\
x&X
\end{array}
\right)\succeq0,\label{SDP}
\end{eqnarray}
where $Q_i\cdot X={\rm trace}(Q_iX)$ is the inner product of matrices $Q_i$ and $X$. Note that the Lagrangian dual problem of problem (P) is
\begin{eqnarray*}
\rm (L)~&\max& \tau\\
&\rm s.t.&\left(\begin{array}{cc}Q_0 &\frac{c_0}{2}\\\frac{c_0^T}{2}&-\tau\end{array}\right)-\sum_i^{l}\lambda_i\left(\begin{array}{cc}Q_i &\frac{c_i}{2}\\\frac{c_i^T}{2}&d_i\end{array}\right)-\sum_i^{m}\mu_j\left(\begin{array}{cc}0 &\frac{a_j}{2}\\\frac{a_j^T}{2}&-b_j\end{array}\right)\succeq0,\\
&&\lambda_i\geq0,~i=1,\ldots,l,~ \mu_j\geq0,~j=1,\ldots,m,
\end{eqnarray*}
which is also known as the Shor's relaxation \cite{shor1987quadratic}. It is well known (see, e.g., \cite{boyd1997semidefinite}) that (L)  is the conic dual  of   (SDP) and (SDP) and (L) have the same optimal value when the strong duality holds for (SDP).  Furthermore, the strong duality holds for  (SDP) when  (SDP) is bounded from below and  Slater condition holds for (SDP). When the Slater condition  holds true for problem (P), i.e., there exists a strictly feasible solution $\hat x$ such that $\hat x^TQ_i\hat x+c_i^T\hat x+d_i<0,~i=1,\ldots,l$ and $a_j^T\hat x\leq b_j,~j=1,\ldots,m$, the Slater condition for (SDP)  automatically holds, e.g., by letting $\hat X=\hat x\hat x^T+\epsilon I$, for sufficiently small $\epsilon>0$  such that   $Q_i\cdot \hat X+c_i^T\hat x+d_i\leq x^TQ_i\hat x+c_i^T\hat x+d_i+\epsilon \lambda_{max}(Q_i) <0$,
where $\lambda_{max}(Q_i)$ is the maximum eigenvalue of matrix $Q_i$.

As the basic SDP relaxation is often too loose, valid inequalities have been considered to strengthen $\rm(SDP)$ in the literature.
One widely used technique in strengthening the basic SDP\ relaxation is  the RLT \cite{sherali2013reformulation}, which  linearizes
the product of any pair of  linear constraints,
i.e.,$$(b_i-a_i^Tx)(b_j-a_j^Tx)=b_{i}b_j-(b_ja_i^T+b_ia_j^T)x+a_i^Txx^{T}a_j\geq0. $$
By linearizing $xx^T$ to $X$, we get a tighter (SDP) relaxation enhanced with the  RLT constraints for problem (P):
\begin{eqnarray}
 {\rm (SDP_{RLT})}
&~ \min&  Q_0\cdot X+c_0^Tx \notag\\
&~\text{s.t.}~&(\ref{QC}),(\ref{LC}),(\ref{SDP}),\notag\\
&~&a_ia_j^T\cdot X+b_ib_j-b_ja_i^Tx-b_ia_j^Tx\geq0,~\forall 1\leq i<j\leq m.\label{RLT}
\end{eqnarray}
Note that when $i=j$, the RLT constraint $a_ia_j^T\cdot X+b_ib_j-b_ja_i^Tx-b_ia_j^Tx\geq0$ is dominated by (\ref{SDP}) and can be omitted.

Moreover, it has been shown in
\cite{burer2012milp} and \cite{sturm2003cones} that  SOC-RLT constraints  can be used to strengthen the convex relaxation $\rm (SDP_{RLT})$ for problem (P).
In particular,  decomposing a positive semidefinite matrix $Q_i$ as $Q_i=B_i^TB_i$, $i\in \mathcal{C}$, we can rewrite the convex quadratic constraint in an  SOC form, i.e.,
\begin{eqnarray}
&&\left.\begin{array}{l}
x^TQ_{i}x\leq-d_i-c_i^Tx  \Rightarrow-d_i-c_i^Tx\geq0\\
x^TQ_{i}x\leq-d_i-c_i^Tx
\end{array}\right\}\Rightarrow\notag \\
&&\norm{\left(\begin{array}{c}B_{i}x\\\frac{1}{2}(-d_i-c_i^Tx-1)\end{array}\right)}\leq\frac{1}{2}(-d_i-c_i^Tx+1).\label{csoca}
\end{eqnarray}
Multiplying the linear term $b_j-a_j^Tx\geq0$ to both sides of the above SOC yields the following valid inequality,
$$(b_j-a_j^Tx)\left(\norm{\left(\begin{array}{c}B_{i}x\\ \frac{1}{2}(1+d_i+c_i^Tx)\end{array}\right)}\right)\leq\frac{1}{2}(b_j-a_j^Tx)(1-d_i-c_i^Tx),$$
whose linearization becomes  the following SOC-RLT constraint,
\begin{eqnarray}
\begin{array}{ll}
&\left\lVert\begin{pmatrix}
B_i(b_jx-Xa_j)\\
\frac{1}{2}(-c_i^TXa_j+(b_jc_i^T-d_ia_j^T-a_j^T)x+(1+d_i)b_j)\end{pmatrix}\right\rVert\\
\leq&\frac{1}{2}(c_i^TXa_j+(d_ia_j^T-a_j^T-b_jc_i^T)x+(1-d_i)b_j),~i\in\mathcal{C},~j=1,\ldots,m.\end{array}\label{CSOC}
\end{eqnarray}

So enhancing  ${\rm (SDP_{RLT})}$ with the SOC-RLT constraints, we get a tighter relaxation for problem ${\rm(P)}$:\begin{eqnarray}
{\rm (SDP_{SOC\mbox{-}RLT})}&\min ~&Q_0\cdot X+c_0^Tx \notag\\
&\text{s.t.} &(\ref{QC}),(\ref{LC}),(\ref{SDP}),(\ref{RLT}),(\ref{CSOC}).\notag
\end{eqnarray}

We have the following theorem due to the obvious inclusion relationship of the feasible regions of the three different relaxations,
${\rm(SDP_{SOC\mbox{-}RLT})},~{\rm(SDP_{RLT})}$ and ${\rm(SDP)}$.
\begin{thm}
~$v{\rm(P)}\geq v{\rm(SDP_{SOC\mbox{-}RLT})}\geq v{\rm(SDP_{RLT})}\geq v{\rm(SDP)}$.
\end{thm}

\subsection{GSRT constraints}

Stimulated by the construction of SOC-RLT constraints, which is only applicable to convex quadratic constraints, we derive the GSRT constraints in this section for general (nonconvex) quadratic constraints.

\subsubsection{GSRT-A constraints}

To construct the GSRT-A constraints for nonconvex quadratic constraints, we first decompose each indefinite  matrix in quadratic constraints according to
the signs of its eigenvalues, i.e., $Q_i=L_i^TL_i-M_i^TM_i$, $i\in\mathcal{N}$, where $L_i$ is  corresponding to the positive eigenvalues and $M_i$ is  corresponding to
the negative eigenvalues. One of such decompositions is the spectral decomposition, $Q_i=\sum_{j=1}^{n-p+r}\lambda_{i_j} v_{i_j}v_{i_j}^T$, where
$\lambda_{i_1}\geq\lambda_{i_2}\cdots\lambda_{i_r}>0>\lambda_{i_{p+1}}\geq \cdots\geq\lambda_{i_n},~0\leq r\leq p<n,,$ and
correspondingly $L_i=(\sqrt{\lambda_{i_1}}v_{i_1},\ldots,\sqrt{\lambda_{i_r}}v_{i_r})^{T},~M_i=(\sqrt{-\lambda_{i_{p+1}}}v_{i_{p+1}},
\sqrt{-\lambda_{i_n}}v_{i_n}).$ A straightforward idea in
applying SOC-RLT is to multiply the linear constraints and the equivalent formula of the nonconvex quadratic constraints resulted from the above decomposition,
\begin{equation}
\norm{\left(\begin{array}{c}L_ix\\\frac{1}{2}(c_i^Tx+d_i+1)\end{array}\right)}\leq \norm{\left(\begin{array}{c}M_ix\\\frac{1}{2}(c_i^Tx+d_i-1)\end{array}\right)},~i\in\mathcal{N}.\label{mlnc}
\end{equation}
Unfortunately, (\ref{mlnc})  is intractable  because of its nonconvexity.
 To overcome this difficulty, we introduce $l-k$ auxiliary variables $z_{i}$, where $l-k$ is the number of nonconvex quadratic constraints, to  replace the right hand side of  (\ref{mlnc}),
\begin{equation*}
z_{i}=\sqrt{x^TM_i^TM_ix+\left(\frac{c_i^Tx+d_i-1}{2}\right)^2}\geq \sqrt{x^TL_i^TL_ix+\left(\frac{c_i^Tx+d_i+1}{2}\right)^2}.
\end{equation*}
We thus get an SOC constraint,
\begin{equation}\label{GSRT-Asaml}\norm{\left(\begin{array}{c}L_ix\\\frac{1}{2}(c_i^Tx+d_i+1)\end{array}\right)}\leq z_{i},
\end{equation} and a nonconvex equality constraint,
\begin{equation}\label{GSRT-Areq}
\norm{\left(\begin{array}{c}M_ix\\\frac{1}{2}(c_i^Tx+d_i-1)\end{array}\right)}=z_{i}.
\end{equation}
We then obtain the following  reformulation of problem $\rm(P)$:
\begin{eqnarray*}
 \rm{(RP)}~~~
&~&\min x^TQ_0x+c_0^Tx\notag \\
&~&\text{s.t.}~x^TQ_ix+c_i^Tx+d_i\leq0,~i=1,\ldots,l,\\
&~&~~~~~\norm{\left(\begin{array}{c}L_ix\\\frac{1}{2}(c_i^Tx+d_i+1)\end{array}\right)}\leq z_{i},~i\in\mathcal{N},\\
&~&~~~~~\norm{\left(\begin{array}{c}M_ix\\\frac{1}{2}(c_i^Tx+d_i-1)\end{array}\right)}=z_{i},~i\in\mathcal{N},\\
&~&~~~~~a_j^Tx\leq b_j,~j=1,\ldots,m.
\end{eqnarray*}

We next construct a convex relaxation by generalizing the SOC-RLT constraints for ${(\rm RP)}$.
First we lift the problem into a matrix space by denoting $\begin{pmatrix}X&S\\S^T&Z\end{pmatrix}=\left(\begin{array}{c}x\\z\end{array}\right)(x^T~z^T)$. We then relax the intractable nonconvex constraint $\begin{pmatrix}X&S\\S^T&Z\end{pmatrix}=\left(\begin{array}{c}x\\z\end{array}\right)(x^{T}\ z^{T})$ to $\begin{pmatrix}X&S\\S^T&Z\end{pmatrix}\succeq\left(\begin{array}{c}x\\z\end{array}\right)(x^{T}\ z^{T})$, which is equivalent
to the following LMI,  by the Schur complement,
$$\left(\begin{array}{ccc}1&x^T&z^T\\x&X&S\\z&S^T&Z\end{array}\right)\succeq0.$$
By multiplying $b_j-a_j^Tx$ and $\norm{L_ix,\frac{1}{2}(c_i^Tx+d_i+1)}\leq z_{i}$, we further get
$$\norm{\left(\begin{array}{c}L_ix(b_j-a_j^Tx)\\\frac{1}{2}(c_i^Tx+d_i+1)(b_j-a_j^Tx)\end{array}\right)}\leq z_{i}(b_j-a_j^Tx),$$
$$\text{i.e.},~\norm{\left(\begin{array}{c}L_ib_jx-L_ixx^Ta_j\\\frac{1}{2}(c_i^T(b_jx-xx^Ta_j)+(d_i+1)(b_j-a_j^Tx))\end{array}\right)}\leq z_{i}b_j-z_{i}x^Ta_j.$$
Then the linearization of the above formula gives rise to
\begin{equation}\label{GSRT-AbigL}
\norm{\left(\begin{array}{c}L_ib_jx-L_iXa_j\\\frac{1}{2}(c_i^T(b_jx-Xa_j)+(d_i+1)(b_j-a_j^Tx))\end{array}\right)}\leq z_{i}b_j-S_{\cdot,i}^Ta_j.
\end{equation}Since  the equality constraint  (\ref{GSRT-Areq}) is nonconvex and intractable, relaxing  (\ref{GSRT-Areq})  to inequality yields the following tractable SOC constraint,\begin{equation}\label{GSRT-Asamr}
\norm{\left(\begin{array}{c}M_ix\\\frac{1}{2}(c_i^Tx+d_i-1)\end{array}\right)}\leq z_{i}.
\end{equation}
 Similarly, we get the following valid inequalities by linearizing the product of
(\ref{GSRT-Asamr}) and $b_j-a_j^Tx$,
\begin{equation}\label{GSRT-AbigR}
\norm{\left(\begin{array}{c}M_ib_jx-M_iXa_j\\\frac{1}{2}(c_i^T(b_jx-Xa_j)+(d_i-1)(b_j-a_j^Tx))\end{array}\right)}\leq z_{i}b_j-S_{\cdot,i}^{T}a_j.
\end{equation}
We also linearize the quadratic form of (\ref{GSRT-Areq}),
\begin{equation*}
\norm{\left(\begin{array}{c}M_ix\\\frac{1}{2}(c_i^Tx+d_i-1)\end{array}\right)}^{2}=z_{i}^2,
\end{equation*}
  to a tractable linearization,\begin{equation}\label{GSRT-Aeq}
Z_{i-k,i-k}=X\cdot M_i^TM_i+\frac{1}{4}(c_ic_i^T\cdot X+(d_i-1)^2+2c_i^Tx(d_i-1)).
\end{equation}
The above constraints connect the variables $Z$, $S$, $X$, $z$ and $x$, which are essential in strengthening the SDP relaxation. Without (\ref{GSRT-Aeq}), $S$, $Z$ and  $z$ would be  unbounded and have no impact on the relaxation.

Finally, (\ref{GSRT-Asaml}), (\ref{GSRT-AbigL}), (\ref{GSRT-Asamr}), (\ref{GSRT-AbigR}) and (\ref{GSRT-Aeq})  together make up the GSRT-A constraints. With the GSRT-A constraint, we  strengthen $\rm(SDP_{RLT})$ to the following tighter relaxation:
\begin{eqnarray*}
{\rm(SDP_{GSRT\mbox{-}A})}
&  \min&Q_0\cdot X+c_0^Tx \notag \\
&~ \text{s.t.} & (\ref{QC}),(\ref{LC}),(\ref{RLT}),(\ref{CSOC}),(\ref{GSRT-Asaml}), (\ref{GSRT-AbigL}), (\ref{GSRT-Asamr}), (\ref{GSRT-AbigR}), (\ref{GSRT-Aeq})  \notag\\
&~&\left(
\begin{array}{ccc}
1&x^T&z^T\\
x&X&S\\
z&S^T&Z\\
\end{array}
\right)\succeq0.
\end{eqnarray*}

The  GSRT-A constraints truly strengthen $\rm (SDP_{SOC-RLT})$  because the projection of the feasible set of problem ${\rm(SDP_{GSRT\mbox{-}A})}$ on ($x$, $X$)
is smaller than the feasible set of $\rm(SDP_{SOC\mbox{-}RLT})$.
From the above paragraph, we know that  GSRT-A constraints consist of five types of constraints: (\ref{GSRT-Asaml}) and (\ref{GSRT-Asamr}) are the new SOC constraints decomposed from the nonconvex quadratic constraints;
 (\ref{GSRT-AbigL}) (respectively,
 (\ref{GSRT-AbigR})) is the linearization of the product of (\ref{GSRT-Asaml}) (respectively, (\ref{GSRT-Asamr}))  and the linear constraints $b_j-a_j^Tx$;
and (\ref{GSRT-Aeq}) is the  linearization of the quadratic form of (\ref{GSRT-Areq}).

The following theorem, which shows the relationship among all the above convex relaxations, is obvious due to the nested inclusion relationship of the feasible regions for this sequence of the relaxations.
\begin{thm}\label{GSRTA}
~$v{\rm(P)}\geq v({\rm SDP_{GSRT\mbox{-}A}})\geq v{\rm(SDP_{SOC\mbox{-}RLT})}\geq v{\rm(SDP_{RLT})}\geq v{\rm(SDP)}$.
\end{thm}

The GSRT-A constraints introduce $2(l-k)\times (m+1)$
extra SOC constraints, where $l-k$ and $m$ are the number of nonconvex quadratic constraints and the number of linear
constraints, respectively, in problem (P),  and the solution process could become time consuming when either or both of $l-k$ and $m$ are large, from which RLT-like methods often suffer. We next present two examples with the same notations as in problem
(P) to show that GSRT-A constraints are possible to achieve a strictly tighter lower bound.

\textbf{Example 1}
$Q_0 =\begin{pmatrix}
0.3 & 0&  0\\
0 & -2   &     0\\
0 &0  &2.4\end{pmatrix}$;  $Q_1= \begin{pmatrix}
1 & 0&  0\\
0 & 1  &     0\\
0 &0  &-1\end{pmatrix}$;
  $a_1=\left(\begin{array}{c}  -0.6\\   -2\\    0.8\end{array}\right)$; $ b_1=-0.5$;
 $c_0=\begin{pmatrix} -0.2  \\  0.8  \\   0.2\end{pmatrix}$; $c_1=0$; $d_1=-1$.

The optimal value is $v{\rm(P)}=-1.21788$ with optimal solution $x^*=(0.05256, 1.00646,-0.125414)^T$.
In this example, $v{\rm(P)}=-1.21788>v{\rm(SDP_{GSRT\mbox{-}A})}=-1.2249>v{\rm(SDP)}=-1.9900 $. A strict inequality holds between $v{\rm(SDP_{GSRT\mbox{-}A})}$
and $v{\rm(SDP)}$.

\textbf{Example 2}
Parameters $Q_0$, $Q_1$, $c_0$, $c_1$, $d_1$, $a_1$ and $b_1$ remain the same as in Example 1, but there is  an extra linear constraint with $a_2=(0.3, 0.2, 0.6)^T$ and $b_2=-0.3$.

The optimal solution is $v{\rm (P)}=-0.7449$ with optimal solution $x^*=( -0.1264,   1.3250,   -0.8785)^T$. In this example, $v{\rm(P)}=-0.7449=v{\rm(SDP_{GSRT\mbox{-}A})}=-0.7449>v{\rm(SDP_{RLT})}=-1.9252>v{\rm(SDP)}=-1.9900$.
A strict inequality holds between ${\rm(SDP_{GSRT\mbox{-}A})}$ and ${\rm(SDP_{RLT})}$.
Moreover, $v{\rm(SDP_{GSRT\mbox{-}A})}=-0.7449$ attains the optimal value, but neither $v{\rm(SDP_{RLT})}=-1.9252$ nor $v{\rm(SDP)}=-1.9900$ does.

Note that the above two examples only involve nonconvex quadratic constraints, so the SOC-RLT constraints are not applicable here. Furthermore,  in Example 1, there are only one linear constraint and one nonconvex quadratic constraint, so the RLT constraints
are not applicable either.

\subsubsection{GSRT-B constraints}

For any type-B constraint satisfying $c_i\in {\rm Range}(Q_i)$, there is an alternative way to express such a nonconvex quadratic constraint,
$$x^TQ_{i}x+c_{i}^Tx+d_{i}=(x+\frac{1}{2}Q_i^{\dagger}c)^TQ_{i}(x+\frac{1}{2}Q_i^{\dagger}c_{i})+d_{i}-\frac{1}{4}c_{i}^TQ_i^{\dagger}c_i.$$
Linearizing the product of the linear term and the SOC constraints generated from type-B nonconvex quadratic constraints
yields the kind of GSRT-B constraints.
Note that this combination fails if  $c_{i}\notin {\rm Range}(Q_i)$,
under which only GSRT\mbox{-}A constraints apply.
For the sake of convenience, we assume type-B constraint holds for all indices $i\in\mathcal{N}, $ in the following of this section.

Using  techniques similar to GSRT\mbox{-}A constraints, we can construct  GSRT-B constraints as follows:
\begin{itemize}
\item i) If $\frac{1}{4}(c_i^TQ_i^{\dagger}c_i)-d_i>0$,  define $\Delta=\sqrt{\frac{1}{4}(c_i^TQ_i^{\dagger}c_i)-d_i}$. We then have the following
type of $\rm GSRT\mbox{-}B$ constraints, termed $\rm GSRT\mbox{-}B_1$ for simplicity,
\begin{eqnarray}
&&~~~~~\norm{L_i(x+\frac{1}{2}Q_i^{\dagger}c_i)}\leq z_i,\label{GSRT-Bl1}\\
&~&~~~~~\norm{\left(\begin{array}{c}M_i(x+\frac{1}{2}Q_i^{\dagger}c_i)\\\Delta\end{array}\right)}\leq z_i,\label{GSRT-Br1}\\
&~&~~~~~Z_{i,i}=M_i^TM_i \cdot (X +\frac{1}{4}Q_i^{\dagger}c_i c_i^TQ_i^{\dagger}+Q_i^{\dagger}c_ix^T)+\Delta^2,\label{GSRT-Beq1}\\
&~&~~~~~\norm{L_i(b_jx-Xa_j+\frac{1}{2}Q_i^{\dagger}c_i(b_j-a_j^Tx))}\leq z_ib_j-a_j^TS_{\cdot,i},\label{GSRT-BbigL1}\\
&~&~~~~~\norm{\left(\begin{array}{c}M_i(b_jx-Xa_j+\frac{1}{2}Q_i^{\dagger}c_i(b_j-a_j^Tx))\\\Delta(b_j-a_j^Tx)\end{array}\right)}\leq z_ib_j-a_j^TS_{\cdot,i},\label{GSRT-BbigR1}\\
&~&~~~~~i\in\mathcal{N},~j=1,\cdots,m;\notag
\end{eqnarray}
\item ii) If $\frac{1}{4}(c_i^TQ_i^{\dagger}c_i)-d_i\leq0$, define $\Delta=\sqrt{d_i-\frac{1}{4}(c_i^TQ_i^{\dagger}c_i)}$.
We then have the following
type of $\rm GSRT\mbox{-}B$ constraints, termed $\rm GSRT\mbox{-}B_2$ for simplicity,
\begin{eqnarray}
&~&~~~~~\norm{\left(\begin{array}{c}L_i(x+\frac{1}{2}Q_i^{\dagger}c_i)\\\Delta\end{array}\right)}\leq z_i,\label{GSRT-Bl2}\\
&~&~~~~~\norm{M_i(x+\frac{1}{2}Q_i^{\dagger}c_i)}\leq z_i,\label{GSRT-Br2}\\
&~&~~~~~Z_{i,i}=M_i^TM_i \cdot (X +\frac{1}{4}Q_i^{\dagger}c_i c_i^TQ_i^{\dagger}+Q_i^{\dagger}c_ix^T),\label{GSRT-Beq2}\\
&~&~~~~~\norm{\left(\begin{array}{c}L_i(b_jx-Xa_j+\frac{1}{2}Q_i^{\dagger}c_i(b_j-a_j^Tx))\\\Delta(b_j-a_j^Tx)\end{array}\right)}\leq z_ib_j-a_j^TS_{\cdot,i},\label{GSRT-BbigL2}\\
&~&~~~~~\norm{M_i(b_jx-Xa_j+\frac{1}{2}Q_i^{\dagger}c_i(b_j-a_j^Tx))}\leq z_ib_j-a_j^TS_{\cdot,i},\label{GSRT-BbigR2}\\
&~&~~~~~i\in\mathcal{N},~j=1,\cdots,m.\notag
\end{eqnarray}
\end{itemize}

For the sake of completeness, we provide a derivation of $\rm(GSRT\mbox{-}B_1)$ as follows:
We first decompose each non-positive   definite matrix in quadratic constraints according to
the signs of its eigenvalues, i.e., $Q_i=L_i^TL_i-M_i^TM_i$, $i\in\mathcal{N}$, as we do for the GSRT\mbox{-}A constraints.
The constraint $x^{T}Q_ix+c_i^Tx+d_i\leq0$ then reduces to
$$(x+\frac{1}{2}Q_i^{\dagger}c_{i})^T(L_i^TL_i-M_i^TM_i)(x+\frac{1}{2}Q_i^{\dagger}c_{i})+d_i-\frac{1}{4}(c_i^TQ_i^{\dagger}c_i)\leq0,$$
and we further have
$$(x+\frac{1}{2}Q_i^{\dagger}c_{i})^T(L_i^TL_i)(x+\frac{1}{2}Q_i^{\dagger}c_{i})\leq (x+\frac{1}{2}Q_i^{\dagger}c_{i})^TM_i^TM_i(x+\frac{1}{2}Q_i^{\dagger}c_{i})+\frac{1}{4}(c_i^TQ_i^{\dagger}c_i)-d_i.$$
Since $\frac{1}{4}(c_i^TQ_i^{\dagger}c_i)-d_i$ is a nonnegative real number and $\Delta=\sqrt{\frac{1}{4}(c_i^TQ_i^{\dagger}c_i)-d_i}$ as defined,
we can then introduce $l-k$ augmented
variables $z_i$ to rewrite the above nonconvex constraints as
\begin{eqnarray*}
z_{i}&=&\sqrt{(x+\frac{1}{2}Q_i^{\dagger}c_{i})^TM_i^TM_i(x+\frac{1}{2}Q_i^{\dagger}c_{i})+\Delta^2}\\
&\geq &\sqrt{(x+\frac{1}{2}Q_i^{\dagger}c_{i})^TL_i^TL_i(x+\frac{1}{2}Q_i^{\dagger}c_{i})},
\end{eqnarray*}
where $l-k$ is the number of nonconvex quadratic constraints. We thus obtain an SOC constraint (\ref{GSRT-Bl1}) from the second inequality, and a nonconvex equality constraint,
\begin{equation}\label{GSRT-Breq1}
\norm{\left(\begin{array}{c}M_i(x+\frac{1}{2}Q_i^{\dagger}c_i)\\\Delta\end{array}\right)}=z_{i}.
\end{equation}
Similarly to the GSRT\mbox{-}A constraints case, we lift the problem by the following matrix inequality,
$$\left(\begin{array}{ccc}1&x^T&z^T\\x&X&S\\z&S^T&Z\end{array}\right)\succeq0.$$
We then obtain (\ref{GSRT-Beq1}) by linearizing the quadratic form of (\ref{GSRT-Breq1}),  i.e.,
$$\norm{\left(\begin{array}{c}M_i(x+\frac{1}{2}Q_i^{\dagger}c_i)\\\Delta\end{array}\right)}^{2}=z_{i}^2.$$
Relaxing the equality in (\ref{GSRT-Breq1}) to inequality yields the SOC constraint (\ref{GSRT-Br1}).
Similar to the GSRT\mbox{-}A\ constraints, by linearizing the product of $b_j-a_j^Tx$ and (\ref{GSRT-Bl1}) ((\ref{GSRT-Br1}), respectively), we further get
the SOC constraint (\ref{GSRT-BbigL1}) ((\ref{GSRT-BbigR1}), respectively).

All the constraints (\ref{GSRT-Bl1}), (\ref{GSRT-Br1}), (\ref{GSRT-Beq1}), (\ref{GSRT-BbigL1}) and (\ref{GSRT-BbigR1}) together make up the $\rm(GSRT\mbox{-}B_1)$ constraints.
The $\rm(GSRT\mbox{-}B_2)$ constraints can be derived in a similar way, whose derivation is omitted for simplicity.

Now we can construct the GSRT-B relaxation for problem $\rm(P)$:
\begin{eqnarray}
{\rm(SDP_{GSRT\mbox{-}B})}~~~
&  \min &Q_0\cdot X+c_0^Tx \notag \\
& \text{s.t.} & (\ref{QC}),(\ref{LC}),(\ref{RLT}),(\ref{CSOC}), \notag\\
&&  \rm(\ref{GSRT-Bl1}-\ref{GSRT-BbigR1})~or~ \rm(\ref{GSRT-Bl2}-\ref{GSRT-BbigR2}),\notag\\
&~& \left(
\begin{array}{ccc}
1&x^T&z^T\\
x&X&S\\
z&S^T&Z\\
\end{array}
\right)\succeq0.\notag
\end{eqnarray}
Similar to Theorem \ref{GSRTA}, the following theorem shows the dominance relationship among different relaxations.
\begin{thm}
~$v{\rm(P)}\geq v({\rm SDP_{GSRT\mbox{-}B}})\geq v{\rm(SDP_{SOC\mbox{-}RLT})}\geq v{\rm(SDP_{RLT})}\geq v{\rm(SDP)}$.
\end{thm}

\begin{rem}
~Although we cannot prove the dominance between GSRT\mbox{-}A and GSRT-B constraints, our numerical experiments show an interesting result: the SDP relaxation enhanced with GSRT-B constraints is always tighter (and faster in most cases) than that enhanced with GSRT\mbox{-}A constraints,
i.e., $v{\rm{(SDP_{GSRT\mbox{-}B})}}\geq v\rm{(SDP_{GSRT\mbox{-}A})}$.
 However, the GSRT\mbox{-}A constraints have their advantages over the GSRT-B constraints, as GSRT\mbox{-}A can be applied to any nonconvex quadratic constraint, while GSRT-B is not applicable to the nonconvex quadratic constraints with $c_i \notin \ra(Q_i)$.
\end{rem}

Note that   the $\rm GSRT\mbox-B_2$ constraint corresponding to index $i$ does not need an auxiliary variable in a special case where $\frac{1}{4}(c_i^TQ_i^{\dagger}c_i)-d_i\leq0$, $M_i(x+\frac{1}{2}Q_i^{\dagger}c_i)$ is a scalar and $M_i(x+\frac{1}{2}Q_i^{\dagger}c_i)\geq0$.
In such a case, the corresponding $\rm GSRT\mbox-B_2$ constraint reduces to
\begin{eqnarray*}
&&\norm{\left(\begin{array}{c}L_i(x+\frac{1}{2}Q_i^{\dagger}c_i)\\\Delta\end{array}\right)}\leq M_i(x+\frac{1}{2}Q_i^{\dagger}c_i),\\
&&\norm{\left(\begin{array}{c}L_i(b_jx-Xa_j+\frac{1}{2}Q_i^{\dagger}c_i(b_j-a_j^Tx))\\\Delta(b_j-a_j^Tx)\end{array}\right)}\leq M_i(b_jx-Xa_j+\frac{1}{2}Q_i^{\dagger}c_i(b_j-a_j^Tx)),\\
&&j=1,\cdots,m,
\end{eqnarray*}
where $\Delta=d_i-\frac{1}{4}(c_i^TQ_i^{\dagger}c_i)>0$.
Under the above conditions, the  relaxation $\rm(SDP_{GSRT\mbox-B})$ reduces to an interesting subcase with a zero duality gap, i.e., minimizing a quadratic function subject to an SOC constraint,
$$x_J^Tx_J\leq (a_1+a_2^Tx)^2,$$
where $x_J$ is a subvector of $x$ with index set $J\subseteq\ \{1,2,\ldots,n\}$, and a special linear
 constraint,
$$a_1+a_2^Tx\geq a_3,$$
where $a_1,a_3\in \Re$  with $a_3>0$ and  $a_2\in \Re^n$, or subject to two special parallel linear constraints,
$$a_4\geq a_1+a_2^Tx\geq a_3,$$
where $a_{4}\in \Re$. This result was first proved, to the best of our knowledge, in \cite{jin2013exact}.

The construction scheme for GSRT-B constraints can also be applied to the convex quadratic constraints if the type-B constraint condition holds, i.e., $c_i\in{\rm Range}(Q_i)$.
 For such type-B convex quadratic constraints, we prove in the following theorem that the  SDP relaxation enhanced with type-B SOC-RLT (SOC-RLT-B) constraints achieves the same optimal value as that enhanced with the conventional SOC-RLT in the literature. On the other hand, the SDP relaxation with SOC-RLT-B constraints  demonstrates a
faster computational speed, which was observed in our numerical tests.
\begin{thm}
~Assume $i\in\mathcal{C}$, $c_i\in {\rm Range}(Q_i)$ and $Q_i\succeq0$, and the following SOC-RLT-B constraint,
\begin{equation}
\norm{B_i(b_jx-Xa_j+\frac{1}{2}Q_i^{\dagger}c_i(b_j-a_j^Tx))}\leq\Delta(b_j-a_j^Tx),\label{A}
\end{equation}
is generated from linearizing the product of $b_j-a_j^Tx\geq0$ and
\begin{equation}\label{BSOCRLT}
\norm{B_ix+\frac{1}{2}Q_i^{\dagger}c_i}\leq\Delta,
\end{equation}
where $\Delta=\sqrt{\frac{1}{4}(c_i^TQ_i^{\dagger}c_i)-d_i}$.
Then the (\ref{A}) is equivalent to the SOC-RLT constraint (\ref{CSOC}).
\end{thm}
\proof
Recall that the SOC-RLT constraint is equivalent to
{\begin{eqnarray*}
&&\norm{B_i(b_jx-Xa_j)}^2+\norm{\frac{1}{2}(-c_i^TXa_j+(b_jc_i^T-d_ia_j^T-a_j^T)x+(1+d_i)b_j)}^2\\
&\leq& \norm{\frac{1}{2}(c_i^TXa_j+(d_ia_j^T-a_j^T-b_jc_i^T)x+(1-d_i)b_j)}^2.
\end{eqnarray*}}
Using the following fact,
{\begin{eqnarray*}
&&\norm{\frac{1}{2}(c_i^TXa_j+(d_ia_j^T-a_j^T-b_jc_i^T)x+(1-d_i)b_j)}^2\\
&-&\norm{\frac{1}{2}(-c_i^TXa_j+(b_jc_i^T-d_ia_j^T-a_j^T)x+(1+d_i)b_j)}^2\\
&=&(b_j-a_j^Tx)(c_i^TXa_j+(d_ia_j^T-b_jc_i^T)x-d_ib_j),
\end{eqnarray*}}
 we obtain $\norm{B_i(b_jx-Xa_j)}^2\leq(b_j-a_j^Tx)(c_i^{T}Xa_j+(d_ia_j^T-b_jc_i^T)x-d_ib_j)$.

Similarly, the SOC-RLT-B constraint (\ref{A}) can be proved to be equivalent to
$$\norm{B_i(b_jx-Xa_j)}^2\leq(b_j-a_j^Tx)(c_i^{T}Xa_j+(d_ia_j^T-b_jc_i^T)x-d_ib_j).$$
\endproof

To summarize, we demonstrated in this section how to construct GSRT\mbox{-}A and GSRT-B constraints to strengthen the SDP relaxations for problem $\rm(P)$.
Numerical tests on these two relaxations will be reported in Section 6 to further verify our theoretical results.

\section{Improvement and extension of the decomposition-approximation method}
In this section, we will introduce an artificial  linear valid inequality for problem (P), which was first proposed by Zheng et al. \cite{zheng2011convex}.
We then propose a new relaxation by introducing  RLT, SOC-RLT and GSRT constraints associated with this new  linear valid inequality and show its
dominance over the decomposition-approximation method in \cite{zheng2011convex}. Adopting the setting in \cite{zheng2011convex} in the following of this
section, we consider problem (P)  with nonnegativity constraint $x\geq0$. To simplify the notations, we include the constraint $x\geq0$
implicitly in the linear constraints $b_j^Tx\leq a_j, ~j=1,\ldots,m$.

Zheng et al. \cite{zheng2011convex} proposed a decomposition-approximation method, by constructing valid inequalities using convex quadratic constraints and an artificial linear constraint.
More specifically, they first introduced an artificial inequality, $\alpha_u=\max\{u^Tx\mid x\in\Omega\}>0$, with a chosen
$u\in\Re_{++}^n=\{y\in\Re^n\mid y_i>0,~i=1,\ldots,n\}$, where $\Omega$ is some suitable set that contains the feasible region.
Although the artificial inequality is redundant itself, it is shown in \cite{zheng2011convex} that the following fact, $$\left(\begin{array}{cc}{\rm diag}(u){\rm diag}(x)&{\rm diag}(u)x\\x^T{\rm diag}(u)&\alpha_u\end{array}\right)\succeq0\Leftrightarrow\alpha_u\geq u^Tx,$$
 yields the following valid LMI that can tighten the SDP relaxation for problem (P),
\begin{equation}\label{alphamatrix}
X\preceq\alpha_u{\rm diag}(u)^{-1}{\rm diag}(x).
\end{equation}
Moreover, using the fact,$$0\succeq\left(\begin{array}{cc}-I_n&B_ix\\x^TB_i^T&c_i^Tx+d_i\end{array}\right) \Leftrightarrow\ x^TB_i^TB_ix+c_i^Tx+d_i\leq0,~i\in\mathcal{C},$$
where $B_i$ is a decomposition of the positive semidefinite matrix $Q_i$ with $Q_i=B_i^TB_i$ as given in Section 2,  the authors in \cite{zheng2011convex}
then developed the following LMI using the Hadamard product,
\begin{eqnarray}
 0&\succeq&\left(\begin{array}{cc}-I_n&B_ix\\x^TB_i^T&c_i^Tx+d_i\end{array}\right)\circ
\left(\begin{array}{cc}{\rm diag}(u){\rm diag}(x)&{\rm diag}(u)x\\x^T{\rm diag}(u)&\alpha_u\end{array}\right)\label{hada}\\
&=&\left(\begin{array}{cc}-{\rm diag}(u){\rm diag}(x)&{\rm diag}(u){\rm Diag}(B_ixx^T)\\
({\rm Diag}(B_ixx^T))^T{\rm diag}(u)&\alpha_u(c_i^Tx+d_i)\end{array}\right).\label{r2sdp}
\end{eqnarray}
Linearizing (\ref{r2sdp}) gives rise to the following HSOC valid inequality,
\begin{eqnarray}
\left(\begin{array}{cc}-{\rm diag}(u){\rm diag}(x)&{\rm diag}(u){\rm Diag}(B_iX)\\
({\rm Diag}(B_iX))^T{\rm diag}(u)&\alpha_u(c_i^Tx+d_i)\end{array}\right)\preceq0.\label{r2sdpl}
\end{eqnarray}
The authors in \cite{zheng2011convex} demonstrated that both constraints in (\ref{alphamatrix})
and (\ref{r2sdpl}) can be used to reduce the relaxation gap of $\rm(SDP_{SOC\mbox-RLT})$. In the following of this section, we will demonstrate that  (\ref{alphamatrix})
and (\ref{r2sdpl}) are redundant for the SDP+RLT+(SOC-RLT) relaxation if we include   $\alpha_u\geq u^Tx$ as an extra linear constraint in problem (P).

We first demonstrate that (\ref{alphamatrix}) is redundant when having RLT constraints associated with $\alpha_u\geq u^Tx$ as an extra linear constraint.
\begin{thm}\label{RLTredun}
~The valid inequality  (\ref{alphamatrix})  is dominated by the RLT constraints generated by $x\geq0$ and $\alpha_{u}\geq u^Tx$,
i.e., $\alpha_ux_i\geq u^TX_{\cdot i}$, $i=1,\ldots,n$.
\end{thm}
\proof
From the RLT constraints derived from $\alpha_{u}\geq u^Tx$ and $x_i\geq0$, i.e, $\alpha_ux_i\geq u^TX_{\cdot i}$, we can conclude
\begin{eqnarray*}
\alpha_u{\rm diag}(u)^{-1}{\rm diag}(x)&=&\left(\begin{array}{ccc}\alpha_u x_1/u_1&&\\&\ddots&\\&&\alpha_u x_n/u_n\end{array}\right)\\
&\succeq&\left(\begin{array}{ccc}u^TX_{\cdot 1}/u_1&&\\&\ddots&\\&&u^TX_{\cdot n}/u_n\end{array}\right).
\end{eqnarray*}
By noting
$$\left(\begin{array}{cc}u_iX_{ij}/u_j&\\&u_jX_{ij}/u_i \end{array}\right)\succeq\left(\begin{array}{cc}&X_{ij}\\X_{ij}& \end{array}\right),~\forall ~1\leq i<j\leq n,$$ and $u^TX_{\cdot,j}=\sum_{i=1}^n u_iX_{ij}$, we immediately have
\begin{eqnarray*}
\alpha_u{\rm diag}(u)^{-1}{\rm diag}(x)\succeq\left(\begin{array}{ccc}u^TX_{\cdot 1}/u_1&&\\&\ddots&\\&&u^TX_{\cdot n}/u_n\end{array}\right)\succeq X,
\end{eqnarray*}
which is exactly (\ref{alphamatrix}).
\endproof

Next, we demonstrate in the following theorem that the HSOC (\ref{r2sdpl}) is redundant when having SOC-RLT constraints.

\begin{thm}\label{hadaredun}
~The HSOC valid inequality (\ref{r2sdpl}) is dominated by the SOC-RLT constraints generated by $x\geq0$,  $\alpha_{u}\geq u^Tx$ and $\norm{B_ix}^2\leq-c_i^Tx-d_i$, i.e., \begin{equation}\norm{\left(\begin{array}{c}
B_iX_{\cdot, j}\\
\frac{1}{2}(x_j+c_i^TX_{\cdot,j}+d_ix_{j})
\end{array}\right)}\leq \frac{1}{2}(x_j-c_i^TX_{\cdot,j}-d_ix_{j})\label{domsoc1}
\end{equation} and
\begin{eqnarray}
&&\norm{\left(\begin{array}{c}\alpha_{u}B_ix-B_iXu\\\frac{1}{2}(\alpha_{u}(1+c_i^Tx+d_i)-(1+d_{i})u^Tx-u^{T}Xc_i)\end{array}\right)}\notag\\
&\leq&
\frac{1}{2}(\alpha_{u}(1-c_i^Tx-d_i)-(1-d_{i})u^Tx+u^{T}Xc_i).\label{domsoc2}
\end{eqnarray}
\end{thm}
\proof
By defining $\frac{0}{0}=0$, due to the Schur complement, (\ref{r2sdpl}) is equivalent to
\begin{equation}\label{hadamard}-\alpha_{u}(c_i^Tx+d_i)\geq \sum\limits_{j=1}^n\frac{(u_jB_{ij}X_{\cdot,j})^2}{u_jx_j},
\end{equation}
where $B_{ij}$ is the $j$th row of the matrix $B_i$. Since $x\geq0$ and $-(c_i^Tx+d_i)\geq x^TB_i^TB_ix=\norm{B_ix}^2$,
we have the SOC-RLT constraints (\ref{domsoc1}),
which is equivalent to
\begin{eqnarray}
\frac{\norm{B_iX_{\cdot,j}}^{2}}{x_{j}}\leq -(c_i^TX_{\cdot,j}+d_ix_{j}).\label{tightsoc}
\end{eqnarray}
From $u>0$, we further have
\begin{eqnarray*}
\frac{u_{j}^2\norm{B_iX_{\cdot,j}}^{2}}{u_{j}x_j}\leq -u_{j}(c_i^TX_{\cdot,j}+d_ix_{j}).
\end{eqnarray*}
Multiplying $u_j$ to both sides of the above inequality and adding the results from 1 to $n$ yield
\begin{eqnarray}
\sum\limits_{j=1}^n\frac{(u_jB_{ij}X_{\cdot,j})^2}{u_jx_j}\leq \sum\limits_{t=1}^n-u_{t}(c_i^TX_{\cdot t}+d_ix_{t})=-(u^{T}Xc_i+d_iu^Tx).\label{general}
\end{eqnarray}
Thus (\ref{general}) implies (\ref{hadamard})  because $-\alpha_{u}(c_i^Tx+d_i)\geq -u^{T}Xc_i+d_iu^Tx$ is hidden
in the SOC-RLT constraint,$$(-\alpha_{u}(c_i^Tx+d_i)+u^{T}Xc_i+d_iu^Tx)(\alpha_u-u^{T}x)\geq\norm{\alpha_{u}B_ix-B_iXu}, $$ which is further equivalent to (\ref{domsoc1}).
We complete our proof by noting the above SOC-RLT constraint is linearized from
$$-\frac{1}{2}(\alpha_{u}-u^Tx)(c_i^Tx+d_i-1)\geq(\alpha_{u}-u^Tx)\norm{B_ix,\frac{1}{2}(c_i^Tx+d_i+1)}.$$
\endproof

In fact,  if the matrix in (\ref{r2sdpl}) is derived from the  SOC
 constraints in any one of (\ref{csoca}), (\ref{BSOCRLT}), (\ref{GSRT-Asaml}), (\ref{GSRT-Asamr}), (\ref{GSRT-Bl1}), (\ref{GSRT-Br1}), (\ref{GSRT-Bl2}) and (\ref{GSRT-Br2}), we can still prove the resulted HSOC valid inequality is redundant.
For simplicity, we term general SOC (GSOC) constraints
for (\ref{csoca}), (\ref{BSOCRLT}), (\ref{GSRT-Asaml}), (\ref{GSRT-Asamr}), (\ref{GSRT-Bl1}), (\ref{GSRT-Br1}), (\ref{GSRT-Bl2}) and (\ref{GSRT-Br2}) and rewrite them in the following unified form,
\begin{equation}\label{GESOC}\norm{C^sx+\xi^s}\leq l_{s}(x,z),~s=1,\ldots,2l-k,
\end{equation}
where $C^s$ can be either $B_i$, $L_i$ or $M_i$ in the above SOC constraints, $\xi^s$ is the corresponding constant  in the norm of the left hand side of the SOC constraints, $l_s(x,z)=(\zeta^s)^Tx+(\eta^s)^Tz+\theta^s$ is a linear function of $x$ and $z$, $\zeta^s\in\Re^n$, $\eta^s\in\Re^{l-k}$ and $\theta^s\in\Re$.
Note that the constraint number $2l-k$ comes from the cardinality of convex constraints, $k$, the number of nonconvex constraints, $l-k$, and the fact that each nonconvex constraint generates two SOC constraints. More specifically, every convex constraint $x^TQ_ix+c_i^Tx+d_i\leq0,~i\in\mathcal{C},$ can be reduced to an SOC constraint in the form of (\ref{GESOC})
with
$l_{i}(x,z)=\frac{1}{2}(-d_i-c_i^Tx+1)$. In particular, we can  set either $l_{i}(x,z)=\frac{1}{2}(-d_i-c_i^Tx+1) $ or $l_{i}(x,z)=1$, if $ c_i\in\ra(Q_i)$.
Besides, every nonconvex constraint $x^TQ_ix+c_i^Tx+d_i\leq0,~i\in\mathcal{N}$, can be reduced to two SOC constraints in the form of (\ref{GESOC})
with $l_{i_1}(x,z)=l_{i_2}(x,z)=z_{i}$ under both type-A or type-B constraint conditions for some $1\leq i_1,i_2\leq 2l-k$.
With a similar analysis, we can extend Theorem \ref{hadaredun} to the following corollary.
\begin{cor}\label{cor1}
~The linearization of the following matrix inequality,
\begin{equation}\label{hadacor}
\left(\begin{array}{cc}l_{s}I&C^sx+\xi^s\\(C^sx+\xi^s)^T&l_{s}\end{array}\right)\circ\left(\begin{array}{cc}{\rm diag}(u){\rm diag}(x)&{\rm diag}(u)x\\x^T{\rm diag}(u)&\alpha_u\end{array}\right)\succeq0,
\end{equation}
is dominated by the GSRT constraints generated by $x\geq0$,  $\alpha_{u}\geq u^Tx$ and $\norm{(C^sx+\xi^s)}\leq l_{s}(x,z)$, $s=1,\ldots,2l-k$.
\end{cor}
\begin{rem}
~In fact, the HSOC valid inequality (\ref{r2sdpl}) is very loose because, from (\ref{tightsoc}), one can find that every entry of the right hand side of (\ref{general}) is
larger than the left hand side of (\ref{general}).\end{rem}

\begin{thm}\label{alphamubound}
~Assume that the relaxation $({\rm SDP_{\alpha GSRT}})$ is obtained by applying RLT, SOC-RLT, and GSRT constraints to problem (P) with a  redundant linear constraint  $u^Tx\leq \alpha_u$.
Then we have $v({\rm SDP_{\alpha GSRT}}) \geq v({\rm SDP_{GSRT}})$ due to the additional valid inequalities in $({\rm SDP_{\alpha GSRT}})$ compared to $({\rm SDP_{GSRT}})$.
\end{thm}
\begin{rem}
~In general, the selected vector $u$ is not necessary to be positive. An interesting research direction is how to identify suitable $u^Tx\leq \alpha_u$ to generate active RLT, SOC-RLT and GSRT constraints.\end{rem}

Next we discuss  two toy examples to show good performance of the relaxation $({\rm SDP_{\alpha GSRT}})$.
The numerical results are shown in Tables 1 and 2. The notation ($\rm SDP$) denotes the basic SDP
 relaxation; ($\rm SDP_{RLT}$) the SDP+RLT relaxation; ($\rm SDP_{SOC\mbox{-}RLT}$) the SDP+RLT+(SOC-RLT) relaxation;
 ($\rm SDP_{\alpha_u}$) ($\rm SDP_{RLT}$) enhanced by (\ref{alphamatrix});  ($\rm SDP_{rtc}$)  ($\rm SDP_{RLT}$)
 enhanced by (\ref{alphamatrix}) and (\ref{r2sdpl}). Moreover, the notation ($\rm SDP_{GSRT\mbox{-}A}$) (($\rm SDP_{GSRT\mbox{-}B}$), respectively) is ($\rm SDP_{SOC\mbox{-}RLT}$) enhanced with GSRT\mbox{-}A constraints (GSRT-B constraints, respectively).
Relaxations (${\rm SDP_{\alpha RLT}}$), (${\rm SDP_{\alpha SOC\mbox-RLT}}$),
(${\rm SDP_{\alpha GSRT\mbox{-}A}}$) and (${\rm SDP_{\alpha GSRT\mbox{-}B}}$)  are (${\rm SDP_{RLT}}$), (${\rm SDP_{SOC\mbox{-}RLT}}$),
(${\rm SDP_{GSRT\mbox{-}A}}$) and (${\rm SDP_{GSRT\mbox{-}B}}$)  enhanced with  RLT,\ SOC-RLT, and GSRT constraints corresponding to  the extra linear constraint  $u^Tx\leq \alpha_u$.
\\
\textbf{Example 3 }\cite{zheng2011convex}
\begin{eqnarray*}
&~  \min& 21x_1^2+34x_1x_2-24x_2^2+2x_1-14x_2 \\
&~  {\rm s.t}&2x_1^2+4x_1x_2+2x_2^2+8x_1+6x_2-9\leq0,\\
&~&-5x_1^2-8x_1x_2-5x_2^2-4x_1+4x_2+4\leq0,\\
&~&x_1+2x_2\leq2,\\
&~&x\in[0,1]^2.
\end{eqnarray*}
\begin{table}[!htbp]
\centering
\caption{SDP bounds for Example 3}
\tabcolsep=2pt
\renewcommand{\arraystretch}{1.3}
\scalebox{0.94}{\begin{tabular}{ll|ll}
\hline \noalign{\smallskip}
SDP relaxation&Lower bound&Extra linear constraint&Lower bound\\
\hline \noalign{\smallskip}
($\rm SDP$)&-20.28&---&---\\
($\rm SDP_{RLT}$)&-16.23&(${\rm \rm SDP_{\alpha RLT}}$)&-11.66\\
($\rm SDP_{SOC\mbox{-}RLT}$)&-13.99&(${\rm SDP_{\alpha SOC\mbox{-}RLT}}$)&-8.445\\
($\rm SDP_{\alpha_u}$)&-10.86&---&---\\
($\rm SDP_{GSRT\mbox{-}A}$)&-6.011&(${\rm SDP_{\alpha GSRT\mbox{-}A}}$)&-4.887\\
($\rm SDP_{GSRT\mbox{-}B}$)& -3.331&(${\rm SDP_{\alpha GSRT\mbox{-}B}}$)& -3.327\\
\hline \noalign{\smallskip}
\end{tabular}}
\end{table}
The optimal value of Example 3 is $v^{*}=-3.327$ with optimal solution $x^*=(0.427,0.588)^T$. In \cite{zheng2011convex}, Zheng et al. set $u=(1,2)^T$, and obtained $\alpha_u=1.8029$.
Strengthening  ($\rm SDP_{SOC\mbox{-}RLT}$) with the decomposition-approximation method, they got
a tighter bound  $v(\rm SDP_{\alpha_u})=-10.86$, compared with ($\rm SDP$), ($\rm SDP_{RLT}$) and ($\rm SDP_{SOC\mbox{-}RLT}$).
We obtain much tighter bounds with our GSRT constraints
 when compared to $(\rm SDP_{\alpha_u})$. The best lower bound -3.327, which is also the optimal value, is achieved by
 (${\rm SDP_{\alpha GSRT\mbox{-}B}}$), i.e., the  combination of RLT, SOC-RLT and GSRT-B constraints with an extra linear constraint $u^Tx\leq \alpha_u$.  It is also remarkable
that ($\rm SDP_{GSRT\mbox{-}B}$) achieves a very good lower bound with -3.331, which demonstrates good performance of GSRT constraints.
\\
\textbf{Example 4} \cite{zheng2011convex}
\begin{eqnarray*}
 &~  \min& -8x_1^2-x_1x_2-13x_2^2-6x_1-x_2 \\
 &~  {\rm s.t}&x_1^2+x_1x_2+2x_2^2-3x_1-3x_2-7\leq0,\\
 &~&2x_1x_2+33x_1+15x_2-10\leq0,\\
 &~&x_1+2x_2\leq6,\\
 &~&x\geq0.
 \end{eqnarray*}
 \begin{table}[!ht]
 \centering
 \caption{SDP bounds for Example 4}
 \tabcolsep=2pt
\renewcommand{\arraystretch}{1.3}
\scalebox{0.94}
  {\begin{tabular}{ll|ll}
\hline \noalign{\smallskip}
 SDP relaxation&Lower bound&Extra linear constraint&Lower bound\\
 \hline \noalign{\smallskip}
 ($\rm SDP$)&-103.43&---&---\\
 ($\rm SDP_{RLT}$)&-26.67& ($\rm SDP_{\alpha RLT}$)&-6.4447\\
 ($\rm SDP_{SOC\mbox{-}RLT}$)&-24.63& ($\rm SDP_{\alpha SOC\mbox{-}RLT}$)&-6.4447\\
 ($\rm SDP_{rtc}$)&-19.61&---&---\\
 ($\rm SDP_{GSRT\mbox{-}A}$)&-24.08& ($\rm SDP_{\alpha GSRT\mbox{-}A}$)&-6.4445\\
 ($\rm SDP_{GSRT\mbox{-}B}$)& -6.4444& ($\rm SDP_{\alpha GSRT\mbox{-}B}$)&-6.4444\\
\hline \noalign{\smallskip}
 \end{tabular}}
 \end{table}
 The optimal value of Example 4 is $v^*=-6.4444$ with optimal solution $x^*=(0,0.6667)^T$.
Zheng et al.  in \cite{zheng2011convex} set $u=(1,1)^T$, obtained $\alpha_u=0.6667$, and got
a tighter bound  $v(\rm SDP_{rtc})=-19.61$, by strengthening  ($\rm SDP_{SOC\mbox{-}RLT}$) with constraints  (\ref{alphamatrix}) and (\ref{r2sdp}).
 For this example, ($\rm SDP_{GSRT\mbox{-}B}$) shows its good quality by achieving a lower bound $-6.4444$ with $x=(0,0.6667)^T$, which is exactly the optimal solution.

The numerical result that ($\rm SDP_{\alpha SOC\mbox{-}RLT}$) is tighter than ($\rm SDP_{\alpha_u}$) and ($\rm SDP_{rtc}$)
verifies the theoretical results in Theorems \ref{RLTredun} and \ref{hadaredun}.
Furthermore, our numerical tests reveal that the GSRT constraints can
improve the quality of the lower bounds when generated with an extra linear constraint $u^Tx\leq \alpha_u$. The fact that in both examples
our relaxations achieve the optimal values demonstrates a good quality of the GSRT\ constraints.


\section{Valid inequalities generated from a pair of SOC constraints}
Recall in Section 2, we construct the GSRT constraint by linearizing the product of an SOC constraint and a linear constraint.
A natural extension is to apply a similar idea to linearize the product of a pair of SOC constraints. However, to the best to our knowledge,
there is no literature that mentions this kind of valid inequalities. In this section, we will  show that valid inequalities generated from the product of
any pair of SOC constraints can indeed tighten the bound for the corresponding SDP relaxation, except for the cases where the two SOC constraints
are both derived from  type-B convex quadratic constraints.

Let us generalize the idea in GSRT constraints to linearize the product of any two SOC constraints.
Multiplying two SOC constraints in the form of (\ref{GESOC}) yields the  valid inequality
\begin{equation}
\norm{C^sxx^{T}(C^t)^T+C^sx(\xi^t)^T+\xi^sx^T (C^t)^T+\xi^s(\xi^t)^T}_F \leq l_{s}l_{t}.\label{socsocrlt}
\end{equation}
Linearizing (\ref{socsocrlt}) yields the following constraint, termed SOC-SOC-RLT (SST) constraint in our paper,
\begin{equation}
\norm{C^sX(C^t)^T+C^sx(\xi^t)^T+\xi^sx^T (C^t)^T+\xi^s(\xi^t)^T}_F \leq \beta_{s,t},\label{sst}
\end{equation}
where $\beta_{s,t}(X,S,Z)=(\zeta^s)^TX\zeta^t+(\zeta^s)^TS\eta^t+(\zeta^t)^TS\eta^s+(\eta^s)^TZ\eta^t+(\theta^s\zeta^t+\theta^t\zeta^s)^Tx+
(\theta^s\eta^t+\theta^t\eta^s)^Tz+\theta^s\theta^t$ is a linear function of  variables $s,z,X,S,Z$, which is linearized from $l_s(x,z)l_t(x,z)$.

Enhanced with valid inequalities (\ref{sst}), we have the following convex relaxation formulation,
\begin{eqnarray*}
\rm{(SDP_{R+SST}) }&~~&\min_{x,X\in R} Q_0\cdot X+c_0^Tx \\
&~&{\rm s.t.}~~{\rm }\norm{C^sX(C^t)^T+C^sx(\xi^t)^T+\xi^sx^T (C^t)^T+\xi^s(\xi^t)^T}_F \leq \beta_{s,t},\\
&~&~~~~~~~~~~~ \forall 1\leq s<t\leq2l-k,
\end{eqnarray*}
where $R$ is the feasible set of either  $\rm(SDP_{GSRT\mbox{-}A})$ or $\rm(SDP_{GSRT\mbox-B})$. Formulation $\rm{(SDP_{R+SST}) }$ introduces $(2l-k)\times(2l-k-1)$
extra matrix norm constraints, which are SOC representable, and thus will be time consuming when  $l$, the number of quadratic constraints, becomes large,
which is a common drawback of RLT-like methods.

The  fact that extra valid inequalities yield a tighter lower bound leads to the following theorem.
\begin{thm}
~$v{\rm(P)}\geq v({\rm SDP_{R+SST})}\geq v(\rm SDP_{R})$.
\end{thm}

To illustrate the SST constraints, consider the following two examples with the same notations in problem (P). For simplicity we only
introduce SST constraints for relaxations with GSRT-A valid inequalities.

\textbf{Example 5 }The parameters in the objective function and quadratic constraints are
$$Q_0=\begin{pmatrix}
41.6520 &   8.7389 &  -3.5465\\
8.7389  &  0.4619 &  13.3579\\
-3.5465 &  13.3579  & 44.4321\\
\end{pmatrix},
Q_1=\begin{pmatrix}
24.2809&    3.5542 &  -5.7609\\
    3.5542&   47.4552&    1.0912\\
   -5.7609 &   1.0912&   36.9438\\
\end{pmatrix},$$
$$Q_2=\begin{pmatrix}
 7.6077 &  16.3267&  -13.0655\\
   16.3267 &  12.6145&  -25.3959\\
  -13.0655 & -25.3959 &   8.0877
\end{pmatrix},Q_3=\begin{pmatrix}
14.3004 &   2.7738 &  12.8803\\
    2.7738  &-18.2473  &  9.5673\\
   12.8803 &   9.5673  &-14.8695
\end{pmatrix},$$
$c_0=\left(\begin{array}{rrr} -45.2696\\   46.8522\\  46.4408\end{array}\right)$,
$c_1=\left(\begin{array}{rrr}  -43.7159\\   23.8375\\   39.8978\end{array}\right)$,
$c_2=\left(\begin{array}{rrr}   -38.1502\\    1.7085\\   37.0175\end{array}\right)$,
$c_3=\left(\begin{array}{rrr}  -31.8133\\  -12.8676\\  -29.7478\end{array}\right)$,\\
$d_1=  -80.4758$, $d_2=   25.4805$, $d_3=   12.1182$, and there is only one linear constraint with
$a=(34.8268,~-22.3518,~-2.6805)^T$, $b=22.0463$.

 Our numerical tests show that for Example 5,
$v(\rm SDP_{GSRT\mbox{-}A})= -21.3379$ and $v(\rm SDP_{GSRT\mbox{-}A+SST})=-21.3151$,
where $(\rm SDP_{GSRT\mbox{-}A})$ is defined in Section 2 and  $(\rm SDP_{GSRT\mbox{-}A+SST})$ is $(\rm SDP_{GSRT\mbox{-}A})$
enhanced with SST constraints (\ref{sst}). Thus, SST constraints indeed tighten the relaxation.

\textbf{Example 6 }The parameters in the objective function and quadratic constraints are
$$Q_0=\left(\begin{array}{rrr}
 21.4825&  -7.7033 &  -0.6240\\
   -7.7033&  -29.8039&   -4.1089\\
   -0.6240 &  -4.1089 &  22.6975\\
\end{array}\right),
Q_1=\left(\begin{array}{rrr}
   37.4987  & -1.0583 &  -1.8307\\
   -1.0583  & 37.1551 &   0.7109\\
   -1.8307  &  0.7109 &  44.4416
\end{array}\right),$$
$$Q_2=\begin{pmatrix}
  -13.5847 &  -0.4516&    4.0519\\
   -0.4516 &  -4.7512&  -17.1011\\
    4.0519 & -17.1011&  -12.0858\\
\end{pmatrix},Q_3=\begin{pmatrix}
 -16.9084 &  18.5030  & 12.8217\\
   18.5030&  -30.1639 &   8.2985\\
   12.8217&    8.2985 & -33.1997
\end{pmatrix},$$
$c_0=\left(\begin{array}{rrr}    34.6975\\    7.5415\\   9.8691\end{array}\right)$,
$c_1=\left(\begin{array}{rrr}   -33.9746\\  -16.6183\\ -23.3710\end{array}\right)$,
$c_2=\left(\begin{array}{rrr}      0.5738\\   41.9009\\  37.4547\end{array}\right)$,
$c_3=\left(\begin{array}{rrr}     40.2865\\   29.6597\\  -44.0517\end{array}\right)$,
$d_1=  -7.0418$, $d_2=  5.4327$, $d_3=    -32.8994$, and there is only one linear constraint with
$a=(    -7.2229,~  45.1322,~  25.0139)^{T}$, $b=37.8832$.

 Numerical tests show that for Example 6,
$v(\rm SDP_{GSRT\mbox{-}A})= -5.51378$ and $v(\rm SDP_{GSRT\mbox{-}A+SST})=-5.3560$.
Thus, SST constraints indeed tighten the relaxation.

The good performance of our relaxation in the above examples demonstrates that the SST constraints can strengthen the SDP relaxation for problem (P) with a significant improvement.

However there is a special case when
the SST\ constraints become being dominated.
In the following we will prove an important theorem to show that (\ref{sst}) is dominated for the basic SDP relaxation
when the two SOC constraints are both derived from two type-B convex
quadratic constraints with $c_i\in \ra(Q_i)$ and $c_j\in \ra(Q_j)$ (where $i$ and $j$ are
the indices of the corresponding convex constraints). This fact could be a main hidden reason why no literature mentions
SST-type valid inequality.
The following lemma helps us prove this result.
\begin{lem}\label{traceinequality}
~If $A$ and $B$ are both $n\times n$ positive semidefinite matrices, then $tr(AB)\leq tr(A)tr(B)$.
\end{lem}
\proof For any vector $u$, let us define $\norm{u}_2=\sqrt{\sum_i u_i^2}$ and   $\norm{u}_1=\sum_i |u_i|$. Since $A$ and $B$ are both positive semidefinite, we have $\norm{\lambda_A}_1=tr(A)$ and  $\norm{\lambda_B}_1=tr(B)$, where $\lambda_A$ and $\lambda_B$ are the vectors formed by all eigenvalues of matrix $A$ and $B$, respectively. We complete the proof using the following fact,
\begin{eqnarray*}
tr(AB)&=&\sum_{i,j} A_{ij}B_{ij}
\leq \norm{A}_F\norm{B}_F\\
&=&\norm{\lambda_A}_2\norm{\lambda_B}_2\leq\norm{\lambda_A}_1\norm{\lambda_B}_1=tr(A)tr(B),
\end{eqnarray*}
where the  first inequality is due to Cauchy-Schwarz inequality.
\endproof

Let us define Type-A SOC constraint if it has the form of (\ref{csoca}), which  can be generated from any convex quadratic constraint,
and Type-B SOC constraint if it has the form of (\ref{BSOCRLT}), which is generated from type-B convex quadratic constraint. Using the above lemma, we will show in the next theorem  that the SST constraints generated by two Type-B SOC constraints that both are derived
from convex constraints are dominated by the linearization of the two associated convex quadratic constraints.
\begin{thm}\label{redundant}
~The SST constraint $$\norm{Q_i^{\frac{1}{2}}XQ_j^{\frac{1}{2}}+Q_i^{\frac{1}{2}}x(\xi^j)^T+\xi^ix^T Q_j^{\frac{1}{2}}+\xi^i(\xi^j)^T}_F\leq l_il_j,$$
 which is generated by $\norm{Q_i^{\frac{1}{2}}(x+\xi^i)}\leq l_i$ and
 $\norm{Q_j^{\frac{1}{2}}(x+\xi^j)}\leq l_j$, $i\neq j$, $i,j\in\mathcal{C}$, is dominated by $$Q_i\cdot X+c_i^Tx+d_i\leq0\text{ and  } Q_j\cdot X+c_j^Tx+d_j\leq0,$$
where $\xi^t=Q_t^{\dagger}c_t$ and $l_t=\frac{1}{4}c_t^TQ_t^{\dagger}c_t-d_t$ is a constant, $t=i~{\rm or}~j$.
\end{thm}
\proof
Define $y=\left(\begin{array}{c}
1\\
x\end{array}\right)$, $Y=\left(
\begin{array}{cc}
1&x^T\\
x&X\end{array}\right)$, $D_i=Q_i^\frac{1}{2}(\xi^i ~I)$ and $D_j=Q_j^\frac{1}{2}(\xi^j ~I)$. Then $\norm{Q_i^{\frac{1}{2}}(x+\xi^i)}\leq l_{i}$ and $\norm{Q_j^{\frac{1}{2}}(x+\xi^j)}\leq l_{j}$ are equivalent to
$\norm{D_iy}\leq l_{i}$ and $\norm{D_jy}\leq l_{j}. $
Also, the SST constraint
$$\norm{Q_i^{\frac{1}{2}}XQ_j^{\frac{1}{2}}+Q_i^{\frac{1}{2}}x(\xi^j)^T+\xi^ix^T Q_j^{\frac{1}{2}}+\xi^i(\xi^j)^T}_F\leq l_{i}l_{j}$$
is equivalent to $\norm{D_iYD_j^T}_F\leq l_{i}l_{j}.$
On the other hand, directly lifting $xx^T$ to $X$ for  $$x^TQ_ix+c_i^Tx+d_i\leq0 \text{ and } x^TQ_jx+c_j^Tx+d_j\leq0$$
yields
$$Q_i\cdot X+c_i^Tx+d_i\leq0\text{ and }Q_j\cdot X+c_j^Tx+d_i\leq0,$$
 which are equivalent to
$tr(D_iYD_i^{T})\leq l_{i}^2$ and $tr(D_jYD_j^{T})\leq l_{j}^2$.

Using the fact
$tr(XY)=tr(YX)$ for any matrix $X\in \Re^{m\times n}$ and $Y\in \Re^{n\times m}$, we complete the proof with the following inequality,
\begin{eqnarray}
\norm{D_iYD_j}_F^2&=&tr((D_iYD_j^{T})^{T}(D_iYD_j^T))\notag\\
&=& tr(D_jY^\frac{1}{2}Y^\frac{1}{2} D_i^TD_iY^\frac{1}{2}Y^\frac{1}{2}D_j^T)\notag\\
&=& tr(Y^\frac{1}{2} D_i^TD_iY^\frac{1}{2}Y^\frac{1}{2}D_j^TD_jY^\frac{1}{2}) \notag\\
&\leq&tr(Y^\frac{1}{2} D_i^TD_iY^\frac{1}{2})tr(Y^\frac{1}{2} D_j^TD_jY^\frac{1}{2})\label{thirdiq}\\
&=&tr(D_iYD_i^T)tr(D_jYD_j^T)\notag\\
&\leq&l_{i}^2l_{j}^2.\notag
\end{eqnarray}
Note that  Lemma \ref{traceinequality} and the fact that $A$ and $B$ are positive semidefinite matrices, where
$A=Y^\frac{1}{2} D_i^TD_iY^\frac{1}{2}$ and $B=Y^\frac{1}{2} D_j^TD_jY^\frac{1}{2}$, are used in the proof of (\ref{thirdiq}).
\endproof

\begin{rem}\label{rem4}
~Note that in Theorem \ref{redundant}, the structure of  $l_t=\frac{1}{4}c_t^TQ_t^{\dagger}c_t-d_i$ indicates that the SOCs are generated from convex quadratic constraints. When the SST valid inequality is generated by two type-A SOC constraints, or a type-A and a type-B SOC constraints and both the SOC constraints are derived from convex constraints,
 our numerical experiments show that the SST valid inequality is still dominated by
 $$Q_i\cdot X+c_i^Tx+d_i\leq0\text{ and  } Q_j\cdot X+c_j^Tx+d_j\leq0,~i,j\in\mathcal{C}.$$
As we are unable to prove the above observation theoretically,  this remains as an open problem in this stage.
\end{rem}

 Note that   in both Examples 5 and 6  the resulted SST constraints are derived from two SOCs at least one of which is not generated from a convex constraint, and our numerical results show that SST constraints indeed help reduce the relaxation gap.
On the other hand, Theorem \ref{redundant} and Remark \ref{rem4} suggest us not to generate SST constraints from two SOCs derived from convex quadratic constraints, in a purpose to avoid generating redundant inequalities.

\section{Valid inequalities in LMI form }
In this section, we introduce and extend  valid inequalities in a form of LMI, i.e., the KSOC valid inequalities, by linearizing the Kronecker products of semidefinite matrices derived from valid SOC constraints, which is motivated by the recent work in  \cite{anstreicher2017kronecker}.
We will further show in this section that these KSOC valid inequalities dominate the HSOC valid inequalities (\ref{r2sdpl}) (which is linearized from (\ref{hada})) and the SST valid inequalities (\ref{sst}) discussed in Sections 3 and 4, respectively.
Moreover, these valid inequalities also shed light on  how to generate valid inequalities that can be easily calculated.

Anstreicher \cite{anstreicher2017kronecker}  introduced a new kind of constraint with an RLT-like technique for the well-known CDT problem
\cite{celis1985trust},\begin{eqnarray*}
&~  \min& x^TBx+b^Tx\notag \\
& ~ \text{s.t.}&\norm{x}\leq1,\\
&~&\norm{Ax+c}\leq 1,
\end{eqnarray*}
where $B$ is an $n\times n$ symmetric matrix and $A$ is an $m\times n$ matrix with full row rank.
By the Schur complement, it is easy to verify that the two quadratic constraints in the CDT\ problem are equivalent to the following LMIs,
\begin{equation}\label{CDT}
\left(\begin{array}{cc}I&x\\x^T&1\end{array}\right)\succeq0~ {\rm and}~
\left(\begin{array}{cc}I&Ax+c\\(Ax+c)^T&1\end{array}\right)\succeq0.
\end{equation}
  Anstreicher \cite{anstreicher2017kronecker} proposed a valid LMI  by linearizing the Kronecker product of the above two matrices, because  the Kronecker product of any two positive semidefinite matrices is positive semidefinite.
To reduce the large dimension of the Kronecker matrix, he further proposed KSOC cuts to handle
the problem of dimensionality.

 We next extend the method in  \cite{anstreicher2017kronecker} to the following two semidefinite matrices,
\begin{equation}\label{geksocm}
\left(\begin{array}{cc}l_s(x,z)I_p&h^{s}(x)\\(h^{s}(x))^T&l_s(x,z)\end{array}\right)~ {\rm and}~
\left(\begin{array}{cc}l_t(x,z)I_q&h^{t}(x)\\(h^{t}(x))^T&l_t(x,z)\end{array}\right),
\end{equation}
 which are derived from (and equivalent to) GSOC constraints in (\ref{GESOC}) by Schur complement, where $h^{j}(x)=C^{j}x+\xi^j$, $j=s,t$.
 We also point out that  the following discussion for (\ref{geksocm}) can also be applied to the case of a pair of two type-A SOC constraints or a type-A SOC constraint and a GSOC constraint, i.e., the following Kronecker product,
$$\left(\begin{array}{cc}-I&B_{i}x\\x^TB_{i}^T&c_i^Tx+d\end{array}\right)\otimes
\left(\begin{array}{cc}l_t(x,z)I_{q}&h^{t}(x)\\(h^{t}(x))^T&l_t(x,z)\end{array}\right).$$
Due to the space consideration, we omit detailed discussion for these cases.

Enlightened by the Kronecker product constraint in  \cite{anstreicher2017kronecker}, we consider the following Tracy--Singh product,
which is just a permutation of the Kronecker product, of the  two matrices in (\ref{geksocm}) (with this reason, we abuse the notation $\otimes$ to denote the Tracy--Singh product for simplicity),
\begin{eqnarray*}
S_{s}&=&\left(\begin{array}{cc}l_{s}(x,z)I_p&h^{s}(x)\\h^{s}(x)^T&l_{s}(x,z)\end{array}\right)\otimes
\left(\begin{array}{cc}l_{t}(x,z)I_q&h^{t}(x)\\h^{t}(x)^T&l_{t}(x,z)\end{array}\right)\notag\\
&=&\left(\begin{array}{cccc}
l_sI_p\otimes l_{t}I_q&l_{s}(x,z)I_p\otimes h^{t}(x)&h^{s}(x)\otimes l_{t}(x,z)I_q& h^{s}(x)\otimes h^{t}(x)\\
 *& l_{s}(x,z)I_p\otimes l_{t}(x,z)&h^{s}(x)\otimes h^{t}(x)^T&h^{s}(x)\otimes l_{t}(x,z)
 \\ *&*&l_{s}(x,z)\otimes l_{t}(x,z)I_q & l_{s}(x,z)\otimes h^{t}(x)\\ *&*&*& l_{s}(x,z)\otimes  l_{t}(x,z)
\end{array}\right),
\end{eqnarray*}
where the notation $*$ is used to simplify the expressions of the entries in the lower triangle which are symmetric to the upper triangle.
Linearizing the above matrix yields
the following KSOC constraint,
\begin{equation}\label{ksoc}\widetilde{S}_{s}=\left(\begin{array}{cccccc}
\beta_{st}I_q&&&K^1&J^{1}&H^1\\
&\ddots&&\vdots&\vdots&\vdots\\
&&\beta_{st}I_q&K^p&J^{p}&H^p\\
*&\cdots&*&\beta_{st}I_p&L^{st}&M^{ts}\\
*&\cdots&*&*&\beta_{st}I_q &M^{st}\\
*&\cdots&*&*&*&\beta_{st}
\end{array}\right)\succeq0,
\end{equation}
where the notations are defined as follows, \\
$M^{st}:=C^tX\zeta^s+C^tS\eta^s+\theta^sC^tx+l_{s}\xi^t$ is a vector linearized from $l_sh^{t}=((\zeta^s)^Tx+(\eta^s)^Tz+\theta^s)(C^tx+\xi^t)$,\\
$K^i:=M^{st}e_i^T$, $i=1,\ldots,p,$ with $e_i$ being the vector with the $i$th entry being 1 and all others being 0s,\\
$J^i:=M^{st}_iI_q$, $i=1,\ldots,p,$\\
$H^i:=C^tX(C_{i,\cdot}^s)^T+\xi_i^sC^tx+C_{i,\cdot}^sx\xi^t+\xi_i^s\xi^t$ is a vector linearized from $h_i^sh^t=(C_{i,\cdot}^sx+\xi_i^s)(C^tx+\xi^t)$,\\
$L^{st}:=C^sX(C^t)^{T}+C^sx(\xi^t)^T+\xi^sC^tx+\xi^s(\xi^t)^T$ is a matrix linearized from $h^{s}(x)\otimes h^{t}(x)^T=h^{s} (h^{t})^T=(C^sx+\xi^s)(C^tx+\xi^t)^T$.

The KSOC cuts in
\cite{anstreicher2017kronecker} remain effective to handle the KSOC constraint  $\widetilde{S}_{s}\succeq0$ when the dimension becomes large.
In addition, an interesting observation is that the SST constraint can be derived from a submatrix of $\widetilde{S}_{s}$.
Specifically, we consider the following submatrix  of $\widetilde{S}_{s}$,
\begin{equation}\label{sstge}
\left(\begin{array}{cccc}
\beta_{st}I_{q}&&&H^1\\
&\ddots&&\vdots\\
&&\beta_{st}I_{q}&H^p\\
*&\cdots&*&\beta_{st}
\end{array}\right).
\end{equation}
By invoking the Schur complement, (\ref{sstge}) yields $\sum_{j=1}^p\frac{{H^j}^TH^j}{\beta_{st}}\leq \beta_{st}.$
With the following fact,
\begin{eqnarray*}&&\sum_{j=1}^p{(H^j})^TH^j=\sum_{j=1}^p\norm{(C^tx+\xi^t)_j(C^sx+\xi^s)}^2\\
&=&\norm{(C^tx+\xi^t)(C^sx+\xi^s)^{T}}_F^2\leq\beta_{st}^2,
\end{eqnarray*}
we conclude that (\ref{sstge}) is equivalent to (\ref{sst}).
Moreover the following matrix inequality,
\begin{equation}\label{subksrt}
\left(\begin{array}{cccc}\beta_{st}I_p&L^{st}&M^{ts}\\
*&\beta_{st}I_q &M^{st}\\
*&*&\beta_{st}
\end{array}\right)\succeq0,
\end{equation}
 which is a submatrix of $\widetilde{S}_{s}$ with a medium size  $(2n+1)\times(2n+1)$, can also be used to tighten relaxations for problem (P).

To summarize, we have invoked the KSOC  constraints in \cite{anstreicher2017kronecker} to derive valid inequalities for
SOC and GSOC constraints. Since the dimension of the Kronecker product matrix increases rapidly as $n$ increases, we intend to adopt computationally cheap valid inequalities via its submatrices to strike a balance between the time cost and  bound quality. More specifically, although (\ref{subksrt}) and SST constraint (\ref{sst}) are submatrices of $\widetilde{S}_{s}$ in (\ref{ksoc}), we may still prefer using these submatrices of KSOC, instead of using (\ref{ksoc}), to generate computationally tractable valid inequalities. We point out that, for a relaxation with a large number of SOC constraints, a practical way is to combine these two methods in an iterative fashion, i.e., solving the relaxation with SST constraints in Section 4  or various submatrices in this section first, and then
finding the Kronecker constraints which violate the semidefiniteness at the current solution $(x,z,X,S,Z),$ and generating
KSOC cuts by the method in  \cite{anstreicher2017kronecker}.



In Section 3, we have demonstrated that the valid inequalities generated by the Hadamard products in (\ref{hada})  and (\ref{hadacor})
are redundant.
In the following, we will generate valid inequalities by replacing the Hadamard products in
 (\ref{hada})  and (\ref{hadacor})
with Kronecker products. Although the  Kronecker product matrices include  the Hadamard product matrices as submatrices (and thus the  corresponding Kronecker product LMIs dominate  (\ref{hada})  and (\ref{hadacor}),
respectively), we will prove that the  two kinds of Kronecker product LMIs  are also redundant.
Let us define
\begin{eqnarray*}
T_{i}&=&\left(\begin{array}{cc}-I&B_{i}x\\x^TB_{i}^T&c_i^Tx+d\end{array}\right)\otimes
\left(\begin{array}{cc}\dg(u)\dg(x)&\dg(u)x\\x^T\dg(u)&\alpha_u\end{array}\right)\\
&=&\left(\begin{array}{cc}-I\otimes \Phi\ &(B_{i}x)\otimes \Phi\\(x^TB_{i}^T)\otimes \Phi&(c_i^Tx+d)\otimes \Phi\end{array}\right),\\
\end{eqnarray*}
where $i\in\mathcal{C}$ and $\Phi=
\left(\begin{array}{cc}\dg(u)\dg(x)&\dg(u)x\\x^T\dg(u)&\alpha_u\end{array}\right)$.
We then define $$V_j^{i}=\left(\begin{array}{cc}\dg(u)\dg(XB_{ij}^T)&\dg(u)XB_{ij}^T)\\B_{ij}X\dg(u)&\alpha_uB_{ij}x\end{array}\right)$$
and$$W^{i}=\left(\begin{array}{cc}\dg(u)\dg(Xc_i+d_i x)&\dg(u)(Xc_{i}+d_{i})\\(Xc_{i}+d_{i})^T\dg(u)&\alpha_u(c_i^Tx+d)\end{array}\right)$$
as linearizations of $(B_{ij}x)\otimes \Phi$ and
$(c_i^Tx+d)\otimes \Phi$, respectively.
Thus linearizing $T_i$ yields the following KSOC valid inequality
\begin{equation}\label{ksoc2}
\widetilde{T}_{i}=\left(\begin{array}{cccc}-\Phi &&& V^i_1\\
&\ddots&&\vdots\\
&&-\Phi& V^i_n\\
*&\cdots&*&W^{i}\end{array}\right)\preceq0.
\end{equation}
One may guess the valid inequality $\widetilde{T}_i\preceq0$ can be used to strengthen relaxations for problem (P) as  $\widetilde{T}_i\preceq0$ dominates the HSOC (\ref{r2sdpl}) (note that (\ref{r2sdpl}) is linearized from  (\ref{hada})), which is a submatrix of $\widetilde{T}_i$.
But, unfortunately, it is redundant, if the relaxation involves SOC-RLT constraints with the artificially introduced redundant linear inequality  $\alpha_{u}\geq u^Tx$, as proved in the following theorem.
\begin{thm}\label{T10}
~The  KSOC inequality $\widetilde{T}_i\preceq0$ is
dominated by the SOC-RLT constraints generated by $x\geq0$,  $\alpha_{u}\geq u^Tx$ and $\norm{B_{i}x}^2\leq -c_i^Tx-d_i$, i.e.,
(\ref{domsoc1}) and (\ref{domsoc2}).
\end{thm}
\proof
Define $P:=\left(\begin{array}{cc}I_{p}&-e\\&1\end{array}\right)$ with $e$ being the all one vector. It is easy  to verify the following facts,
\begin{eqnarray*}
&&\Phi':=P^T\Phi P=\left(\begin{array}{cc}\dg(u)\dg(x)&\\&\alpha_u-u^Tx\end{array}\right),\\
&&{V_j^i}':=P^TV_j^iP=\left(\begin{array}{cc}\dg(u)\dg(XB_{ij}^T)&\\&\alpha_uB_{ij}x-u^{T}XB_{ij}^T\end{array}\right),\\
&&{W^i}':=P^TW^{i}P=\left(\begin{array}{cc}\dg(u)\dg(Xc_i+d_i x)&\\&\alpha_u(c_i^Tx+d)-u^{T}(Xc_i+d_i x)\end{array}\right).
\end{eqnarray*}
So we have the following transformation,
\begin{equation}\label{gksoc}
(I\otimes P)^T\widetilde{T}_{i}(I\otimes P)=\left(\begin{array}{cccc}-\Phi' &&& {V_1^i}'\\
&\ddots&&\vdots\\
&&-\Phi'& {V^i_n}'\\
*&\cdots&*&{W^i}'\end{array}\right).
\end{equation}
From the Schur complement, $(\ref{gksoc})\preceq0$ is equivalent to
${W^i}'\preceq0$ and
\begin{eqnarray*}
\bar{T}_i&:=&{W^i}'-({V_1^i}'~\cdots ~{V_n^i}')\dg(-\Phi',\cdots,-\Phi')^{\dagger}({V_1^i}'~\cdots ~{V_n^i}')^T\\
&=&{W^i}'+\sum_{j=1}^n {V_j^i}'\Phi'^{\dagger}{V_j^i}\preceq0.
\end{eqnarray*}
Together with the fact that $\bar{T}_i$ is a diagonal matrix (since ${W^i}'$, $\Phi'$ and ${V_j^i}'$ are all diagonal), $\bar{T}_i\preceq0$ is equivalent to
$$\alpha_u(c_i^Tx+d)-u^{T}(Xc_i+d_i x)+\frac{\sum_{j=1}^{n}(\alpha_uB_{ij}x-u^{T}XB_{ij}^T)^{2}}{\alpha_u-u^{T}x}\leq0 $$
and
$$u_t(Xc_i+d_ix)_t+\frac{\sum_{j=1}^n [u_t(XB_{ij})_t]^{2}}{u_tx_t}\leq0,~t=1\ldots,n,$$
The former equation is equivalent to (\ref{domsoc2}), and the latter equations are equivalent to, by eliminating $u_t$, (\ref{domsoc1}).
\endproof
Similarly we have the following result for the KSOC constraint generated from   a GSOC and $\Phi$. Although the KSOC constraint dominates the HSOC constraint generated by (\ref{hadacor}), the KSOC constraint is redundant when having GSRT constraints.
\begin{cor}\label{cor2}~The KSOC constraint generated by the following Kronecker product
\begin{equation}\label{Eq54}
\left(\begin{array}{cc}l_{s}(x,z)I&h^s(x)\\(h^{s}(x))^T&l_{s}(x,z)\end{array}\right)\otimes
\left(\begin{array}{cc}\dg(u)\dg(x)&\dg(u)x\\x^T\dg(u)&\alpha_u\end{array}\right)
\end{equation}
is dominated by the GSRT constraints generated by $x\geq0$,  $\alpha_{u}\geq u^Tx$ and $\norm{(C^sx+\xi^s)}\leq l_{s}(x,z)$.
\end{cor}

With a similar analysis, we can prove the KSOC constraint generated by the following Kronecker product \begin{eqnarray}
\left(\begin{array}{cc}l_{s}(x,z)I&h^s(x)\\(h^{s}(x))^T&l_{s}(x,z)\end{array}\right)\otimes
\left(\begin{array}{cc}\dg(u)\dg(x)&\dg(u)x\\x^T\dg(u)&u^{T}x\end{array}\right)\succeq0
\end{eqnarray}
is also dominated by GSRT constraints generated from $x\geq0$, and $\norm{(C^sx+\xi^s)}\leq l_{s}(x,z)$.

In summary, we have demonstrated that the two valid inequalities generated by the Kronecker products in (\ref{ksoc2}) and (\ref{Eq54}) are redundant, although they are more general than
the  associated HSOC constraints (\ref{r2sdpl}) in Theorem \ref{hadaredun} and (\ref{hadacor}) in Corollary \ref{cor1}.

\section{Numerical results}
\begin{figure}[tb]\label{fig1}
\begin{minipage}{0.49\linewidth}
\centering
\includegraphics[width=7cm]{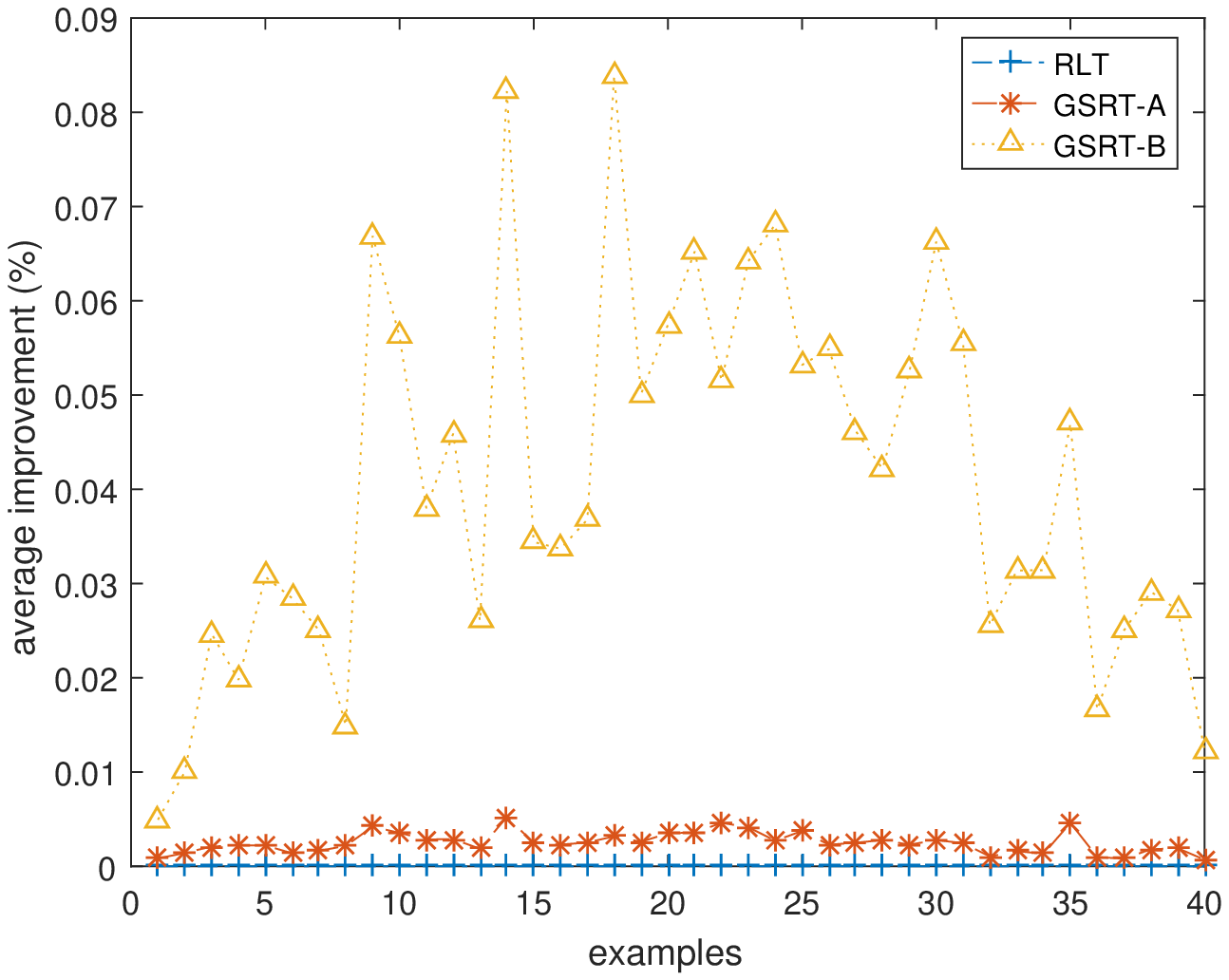}
\end{minipage}%
\hspace*{0.2cm}
\begin{minipage}{0.49\linewidth}
\centering\includegraphics[width=7cm]{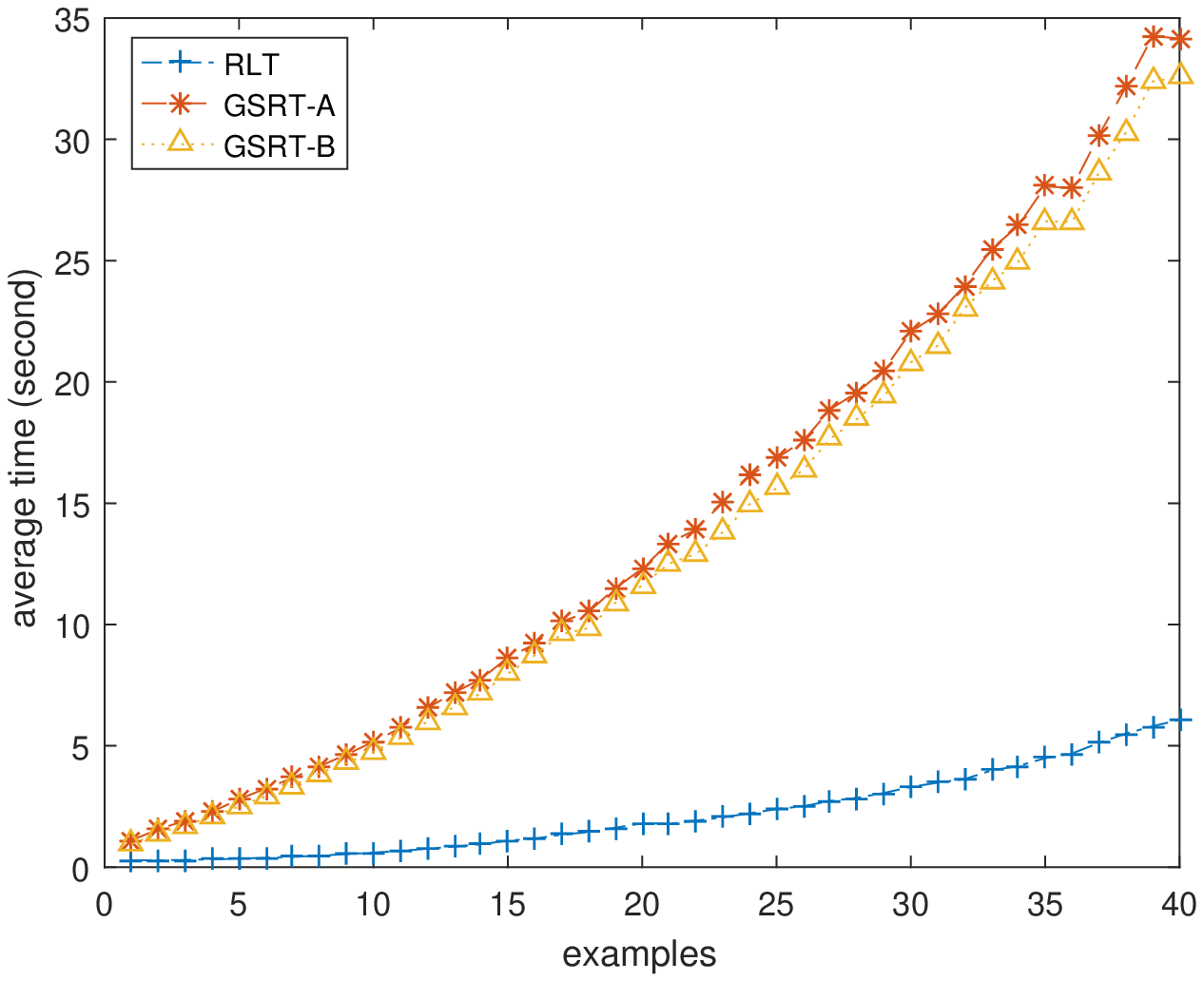}
 \end{minipage}

\begin{minipage}{0.49\linewidth}
\centering
\includegraphics[width=7cm]{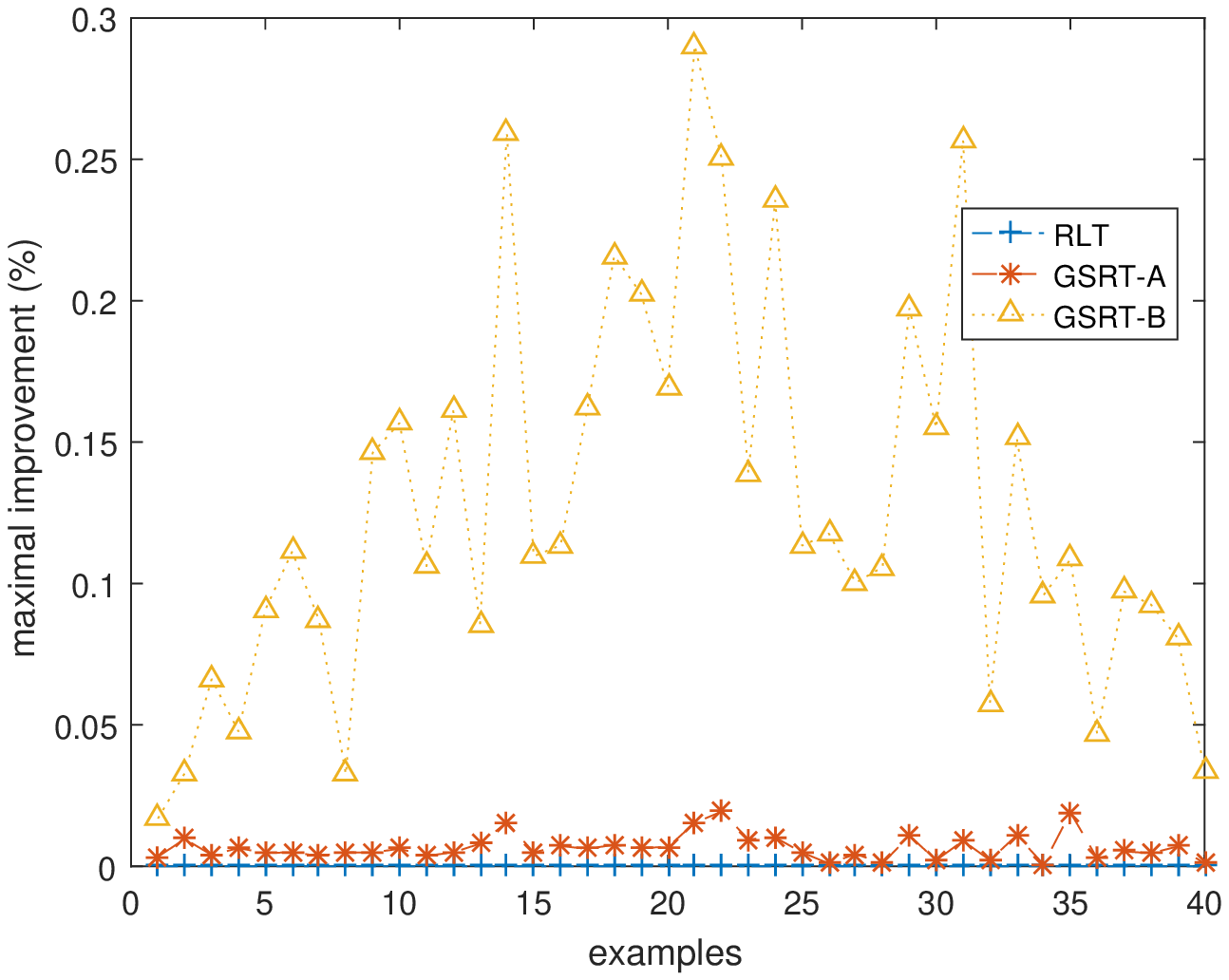}
\end{minipage}%
\hspace*{0.2cm}
\begin{minipage}{0.49\linewidth}
\centering\includegraphics[width=7cm]{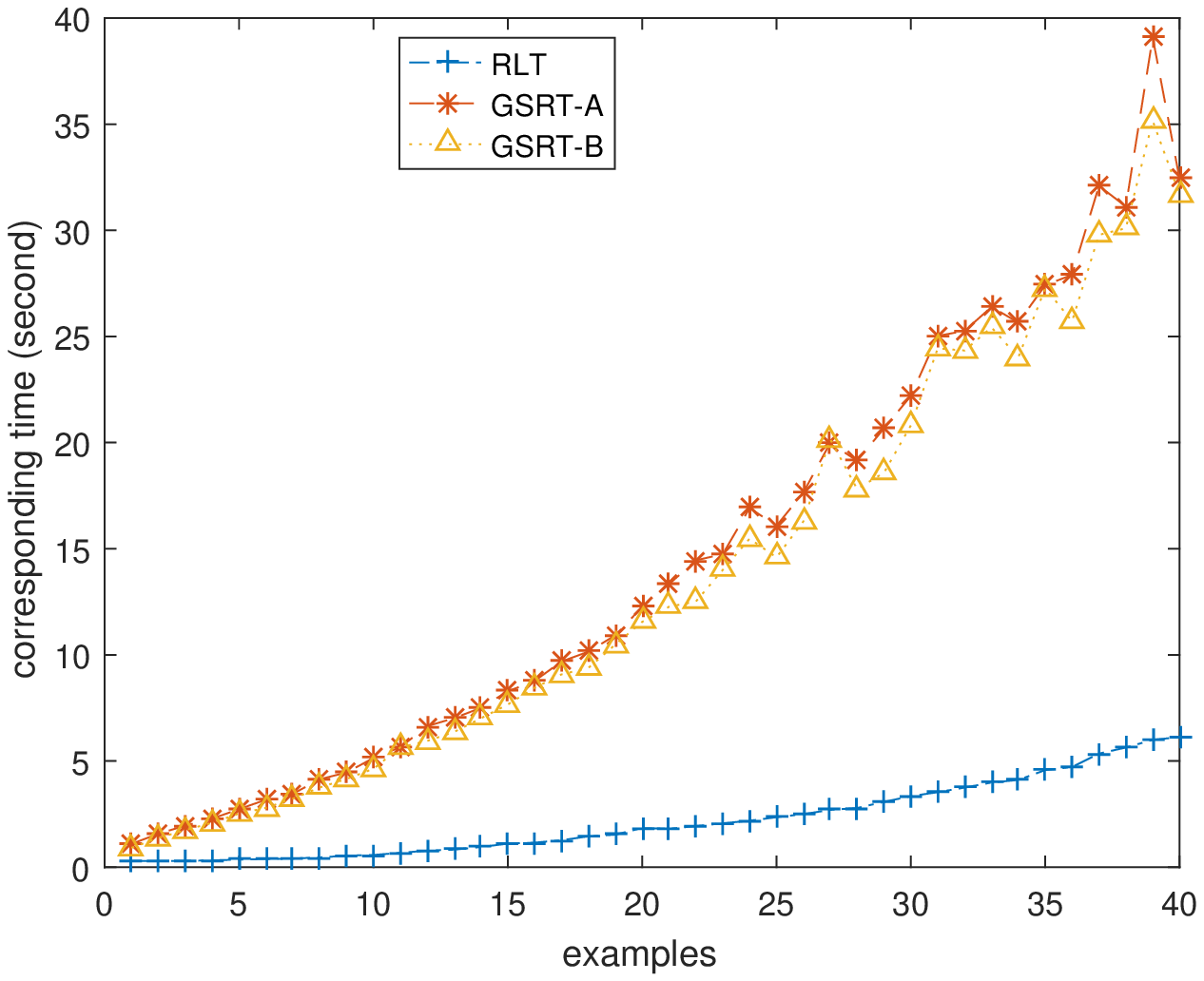}
 \end{minipage}
 \caption{Evolution of average and maximal improvement (of 10 examples) versus number of linear constraints for problem setting $n=20$, $\phi=5$.}\label{figeg5}
\end{figure}

\begin{figure}[htb]\label{fig2}
\begin{minipage}{0.49\linewidth}
\centering
\includegraphics[width=7cm]{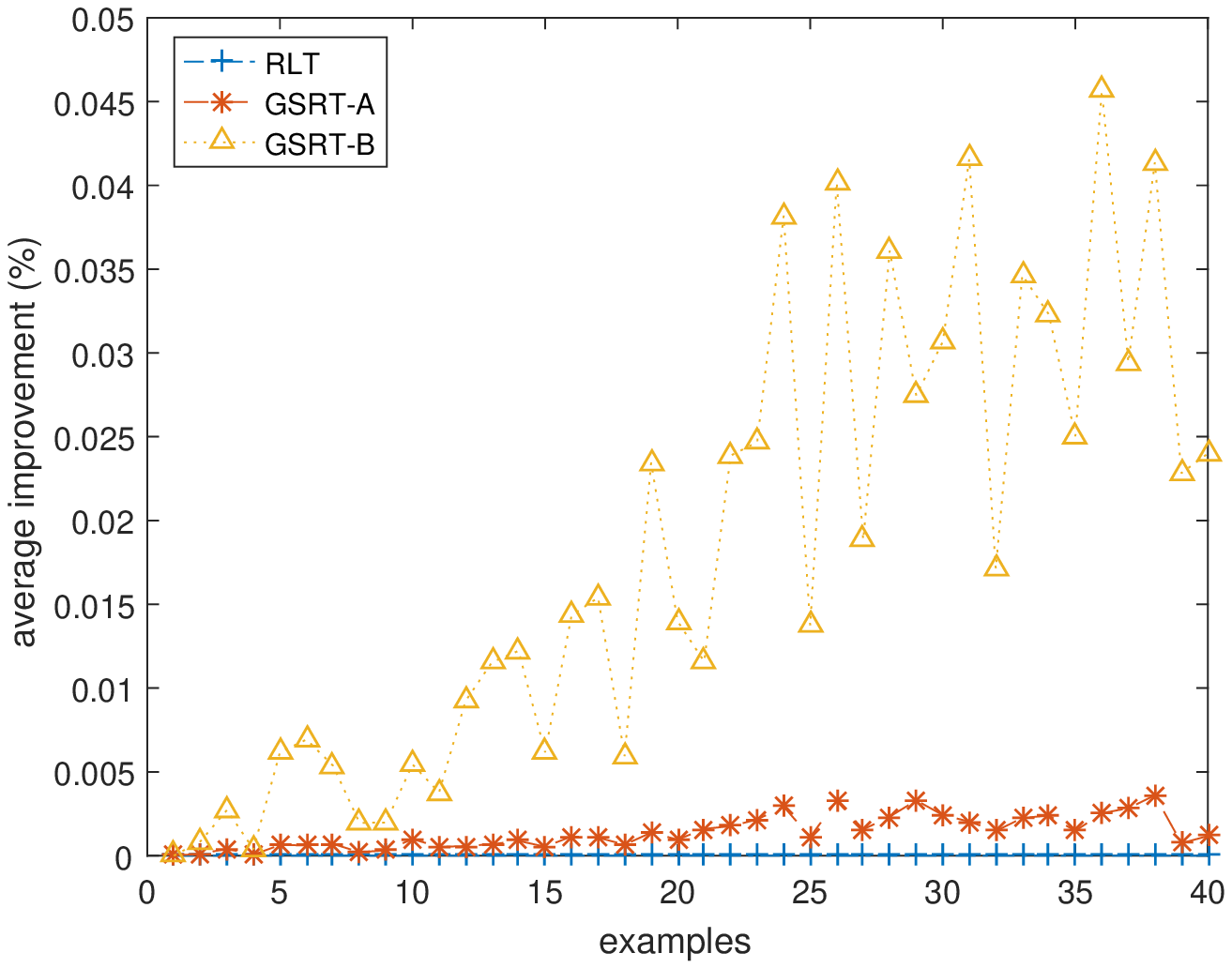}
\end{minipage}%
\hspace*{0.2cm}
\begin{minipage}{0.49\linewidth}
\centering\includegraphics[width=7cm]{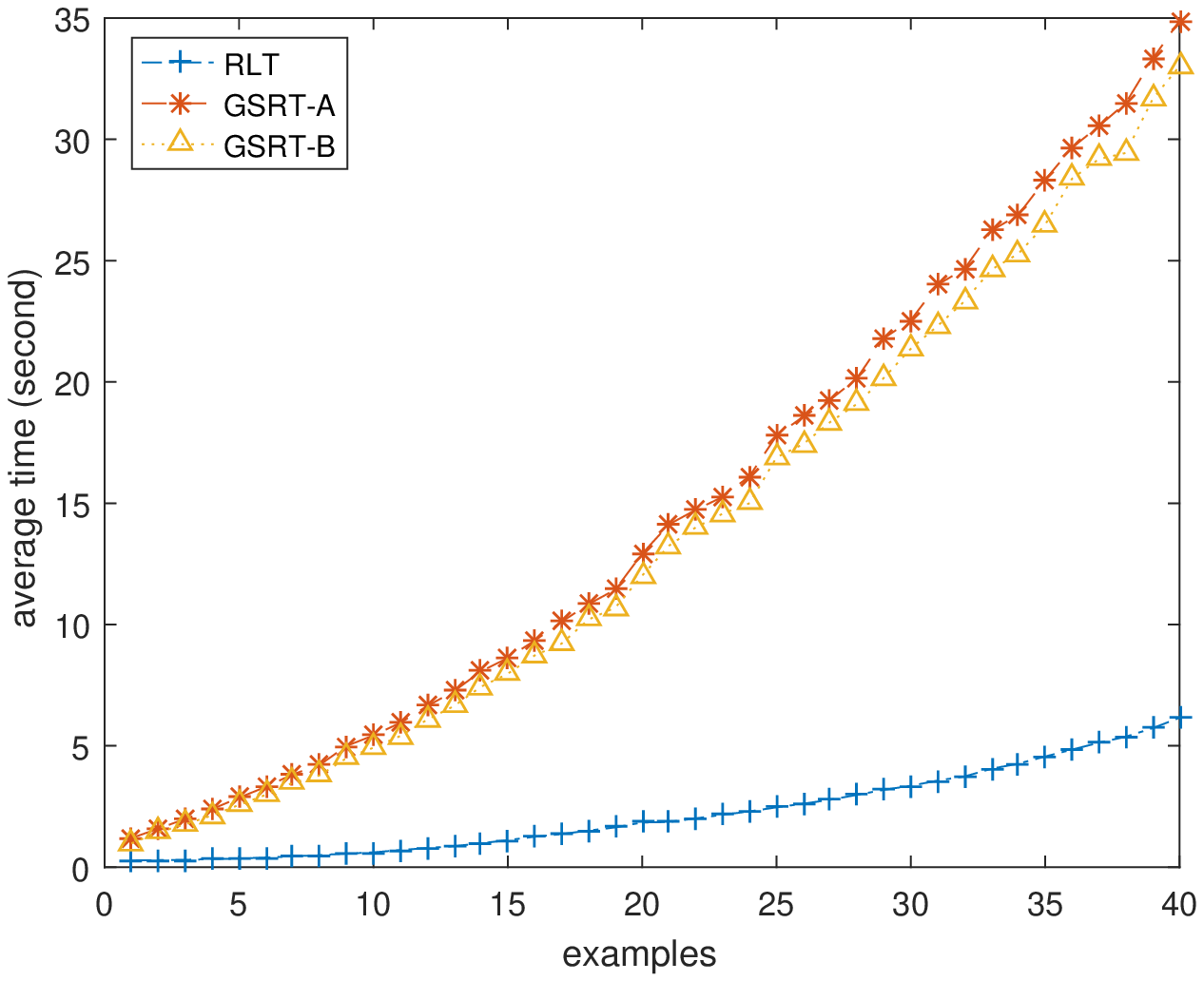}
 \end{minipage}
\begin{minipage}{0.49\linewidth}
\centering
\includegraphics[width=7cm]{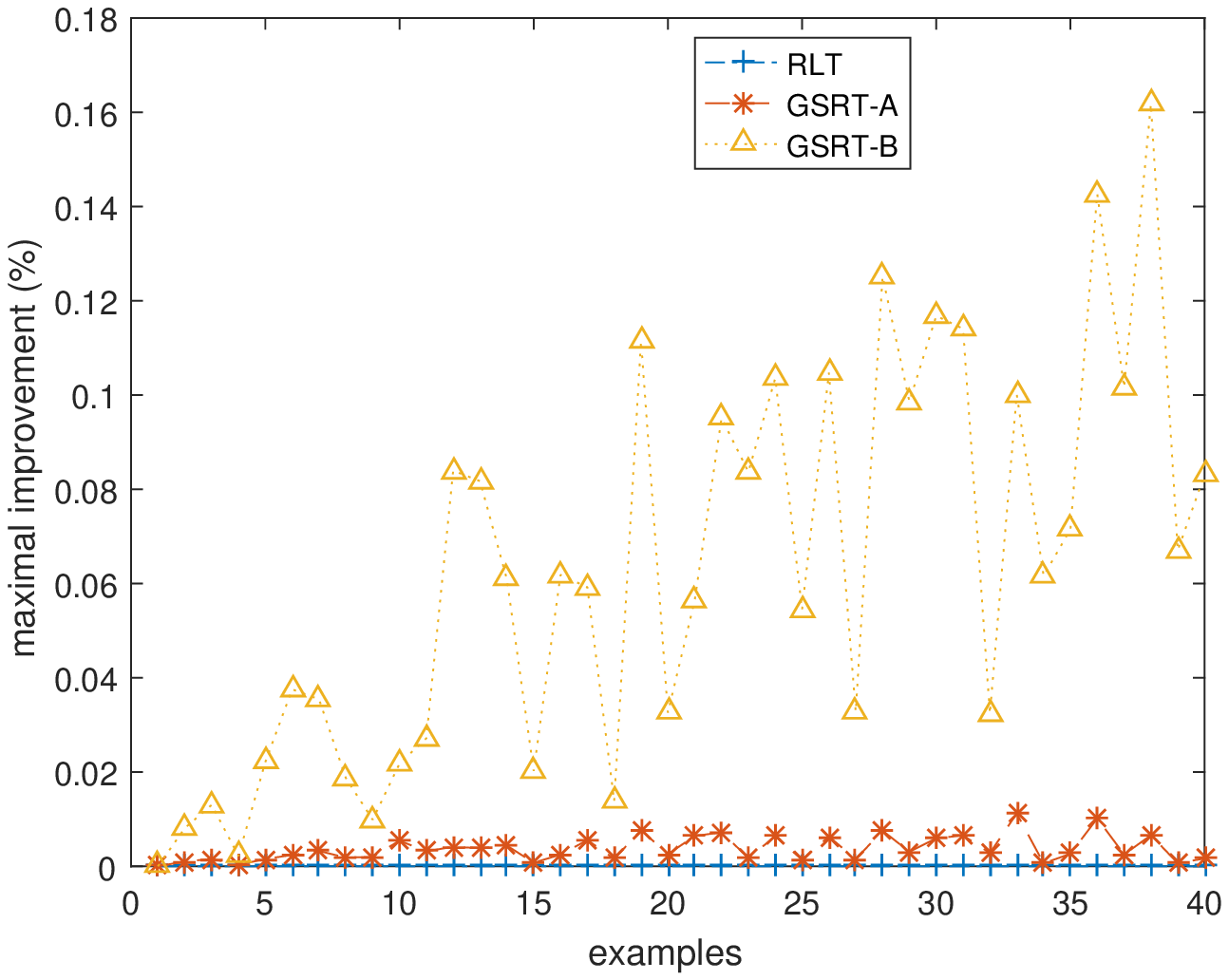}
\end{minipage}%
\hspace*{0.2cm}
\begin{minipage}{0.49\linewidth}
\centering\includegraphics[width=7cm]{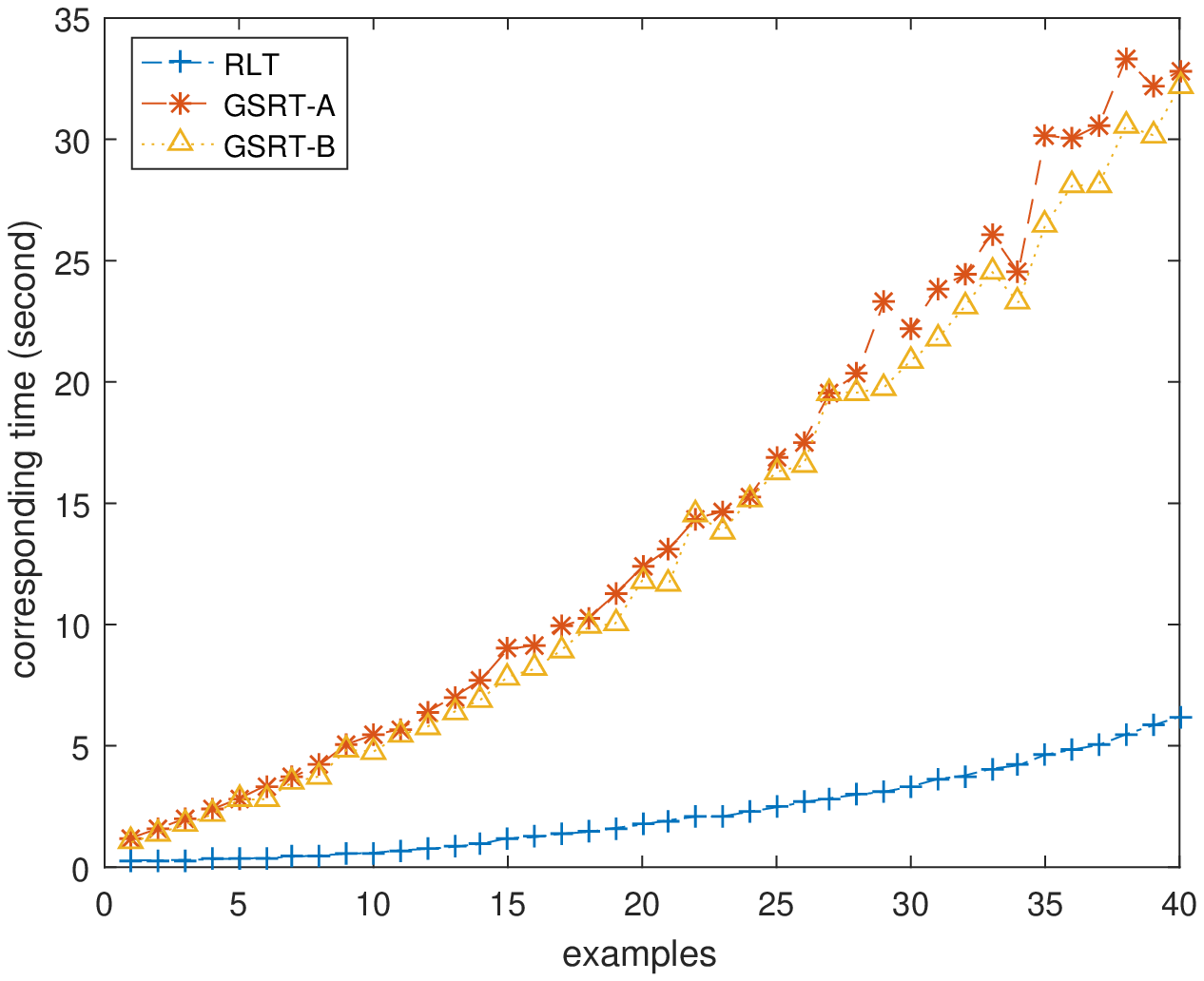}
 \end{minipage}
 \caption{Evolution of average and maximal improvement  versus number of linear constraints for problem setting  $n=20$, $\phi=10$.  }\label{figeg10}
\end{figure}

\begin{figure}[htb]\label{fig2}
\begin{minipage}{0.49\linewidth}
\centering
\includegraphics[width=7cm]{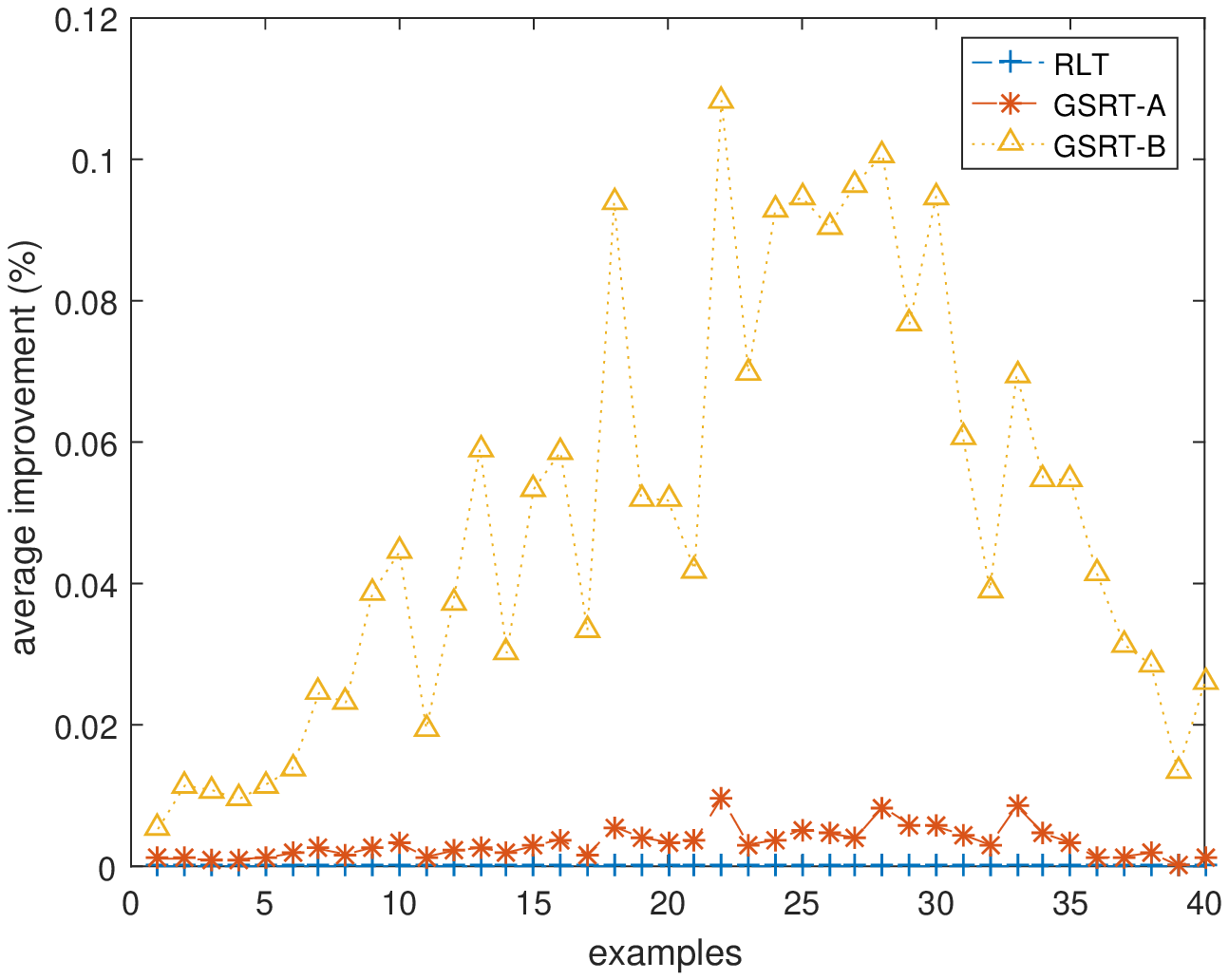}
\end{minipage}%
\hspace*{0.2cm}
\begin{minipage}{0.49\linewidth}
\centering\includegraphics[width=7cm]{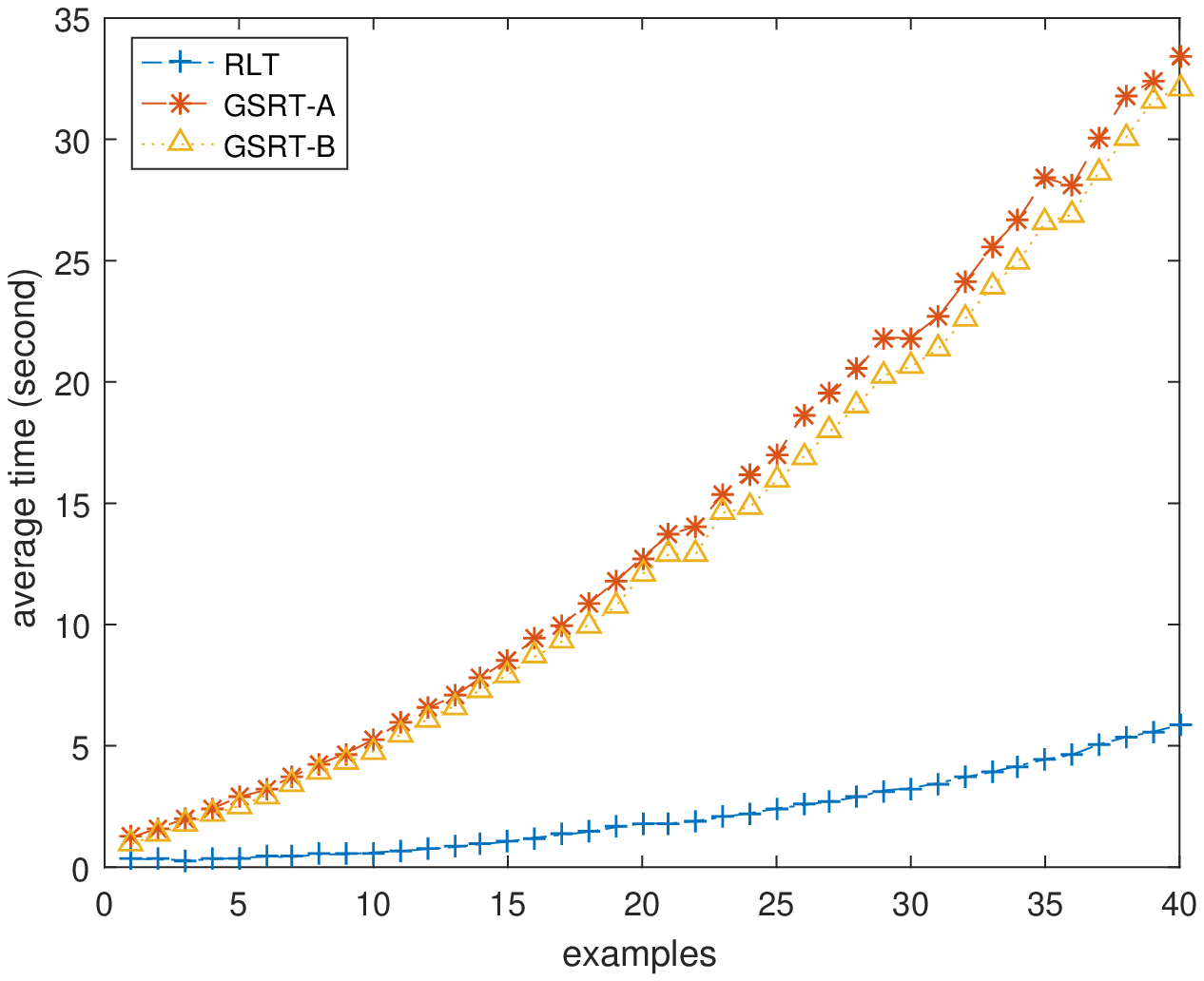}
 \end{minipage}
\begin{minipage}{0.49\linewidth}
\centering
\includegraphics[width=7cm]{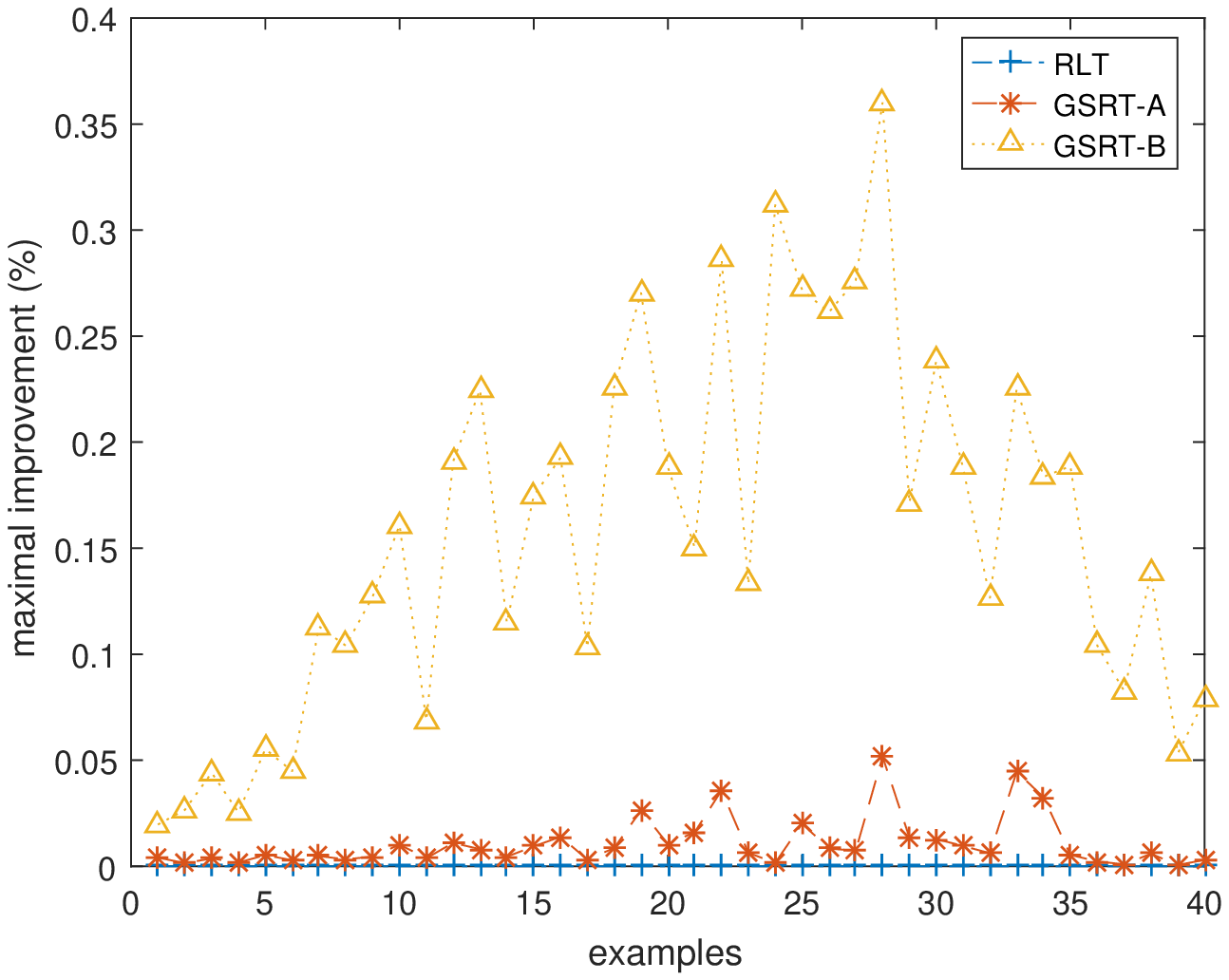}
\end{minipage}%
\hspace*{0.2cm}
\begin{minipage}{0.49\linewidth}
\centering\includegraphics[width=7cm]{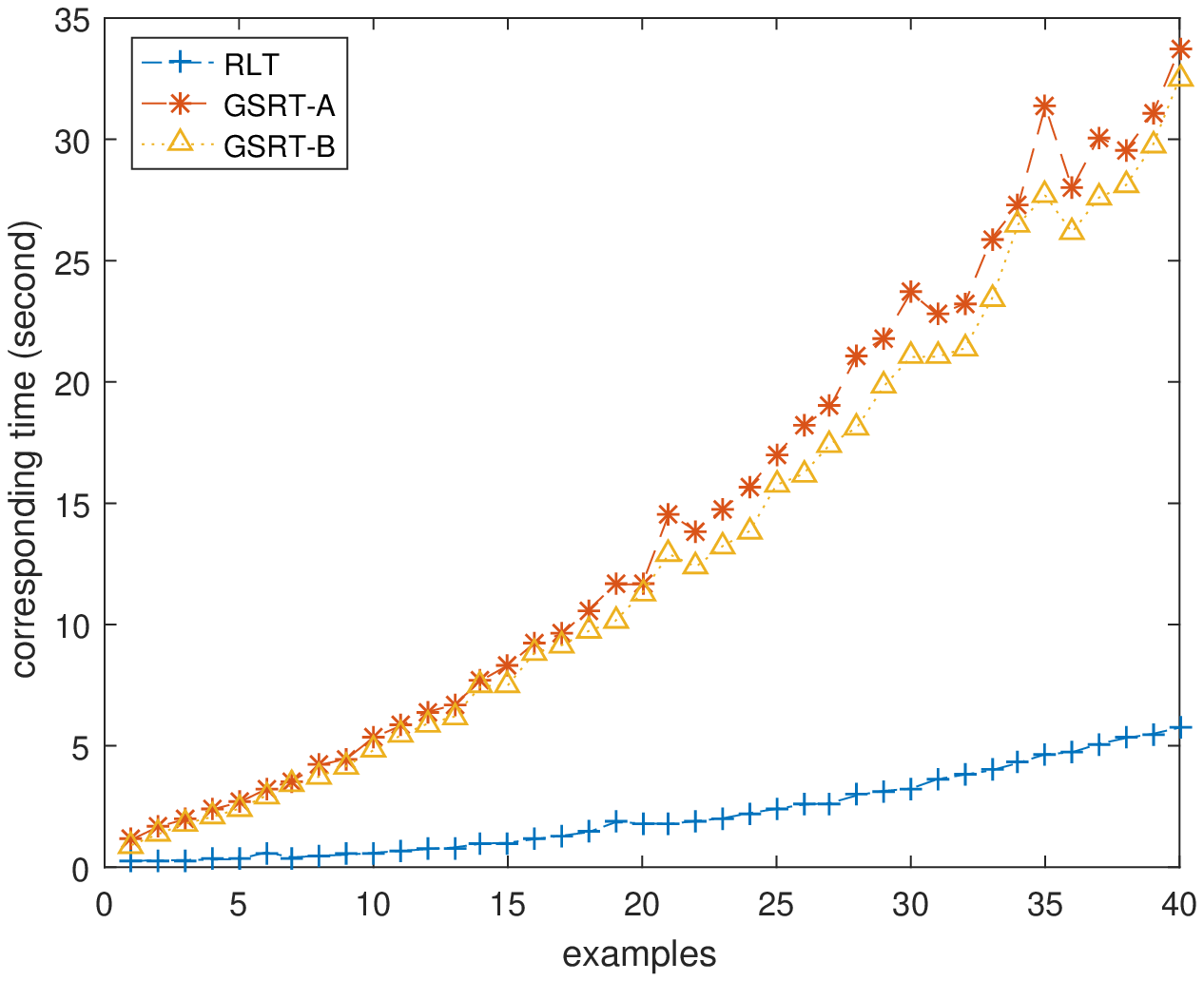}
 \end{minipage}
 \caption{Evolution of average and maximal improvement  versus number of linear constraints for problem setting $n=20$, $\phi=15$.  }\label{figeg15}
\end{figure}

In this section, we report our numerical tests on SDP bounds generated by $\rm{(SDP_{RLT})}$, $\rm{(SDP_{SOC\mbox{-}RLT})}$ and $\rm{(SDP_{GSRT})}$.
The numerical tests in Table \ref{tab3} were implemented in Matlab 2013a, 64bit and was run on a Linux machine with 48GB RAM, 2.60GHz cpu and 64-bit CentOS Release 5.5. And the numerical tests in Figures \ref{figeg5}--\ref{figeg15} were implemented in Matlab2016a  and was run on a PC with 8GB RAM, 3.30GHz cpu and 64-bit Windows 7.
The mixed SDP and SOCP problems in all our numerical examples are modeled by CVX 2.1 \cite{gb08,cvx},  and solved by SDPT3 4.0 within CVX.

 The examples in Table \ref{tab3} were generated  in the following way,  which is similar to Set 1 in \cite{zheng2011convex} but without the box constraint $[0,1]^n$.
The test problems have nonconvex objective function, $k$ convex quadratic constraints, $l-k$ nonconvex quadratic constraints and $m$ linear constraints.
In the following, we use $\xi\in_u [a,b]$ to represent a random number $\xi$ uniformly distributed in the
interval $[a,b]$ and $\rm round(\cdot)$ to represent the value after rounding for a matrix, vector, or scalar. To invoke the GSRT-B valid
inequalities, we choose the instances whose  nonconvex quadratic constraints correspond to   nonsingular matrices.
\begin{itemize}
\item  $Q_0={\rm round }(P_0T_0P_0)$,
$Q_i=P_iT_iP_i$ ($1\leq i\leq l$); $P_i=U_{i1}U_{i2}U_{i3}$,
$U_{it}=I-2\frac{w_tw_t^T}{\|w_t\|^2}$, $i=0,\ldots,l,$ $t=1,2,3$,
$w_t=(w_{t1},\ldots,w_{tn})^T$, $w_{tk}\in_u [-1,1]$.
\item For $1\leq i\leq k$,   $T_i={\rm diag}(T_{i1},\ldots,T_{in})$ with $T_{it}\in_u [0,50]$, for $t=1,\ldots,n$;
For $k+1\leq i\leq l$, $T_{it}\in_u [-50,0]$ for $t=1,\ldots,\frac{n}{2}$ and $T_{it}\in_u [0, 50]$ for
$t=\frac{n}{2}+1,\ldots,n$; $T_{0t}\in_u[-50,50]$, for $t=1,\ldots,n$. Also,
$c_i=(c_{i1},\ldots,c_{in})^T$ with $c_{0t}
 \in_u [-50,50]$, $c_{it}\in_u [-100,0]$ for  $1\leq i\leq k$ and $c_{it}\in_u [0,100]$ for $k+1\leq i\leq l$,  $t=1,\ldots,n$. And $d_i\in_u [-100+\theta^i,\theta^i]$ for  $1\leq i\leq k$ and $d_i\in_u [-10+\theta^i,\theta^i]$ for
 $k+1\leq i\leq l$, where $\theta^i=-e^T_1Q_ie-b_i^Te_1$ with $e_1=(1,0,\ldots,0)^T$.
 \item For $1\leq j\leq m$, $a_j={\rm round}(a_{j1},\ldots,a_{jn})^T$, $a_{jt}
\in_u [-50,50]$, $b_j={\rm round}(\theta_j)$, where $\theta_j\in_u [-10-\vartheta_j,-\vartheta_j]$ with $\vartheta_j=0.5\sum_{j=1}^n\max\{0,a_{jt}\}$, for $t=1,\ldots,n$.
\end{itemize}

\begin{table}[ht!]\small
\centering
\caption{Numerical tests for different convex relaxations}\label{tab3}
\tabcolsep=2pt
\renewcommand{\arraystretch}{1.3}
\scalebox{0.94}
{\begin{tabular}{l|llll|llll}
\toprule\hline  \noalign{\smallskip}
Instance&\multicolumn{4}{l}{Lower bound}&\multicolumn{4}{l}{CPU time}\\\cline{2-5}\cline{6-9}
&RLT&SOC-RLT&GSRT\mbox{-}A&GSRT-B&RLT&SOC-RLT&GSRT\mbox{-}A&GSRT-B\\
\hline  \noalign{\smallskip}
set-30-2-1-59&-972.354&-971.983&-971.836&-971.346&68.5823&110.394&273.571&243.123\\
\hline  \noalign{\smallskip}
set-30-3-1-6&-6049.13&-4650.05&-4635.05&-4497.73&1.73537&5.69738&15.9725&13.4898\\
set-30-3-2-20&-901.782&-890.771&-890.474&-882.626&12.769&28.4496&63.3376&53.7678\\
\hline  \noalign{\smallskip}
set-30-4-1-27&-3697.12&-3574.71&-3573.84&-3497.22&23.4023&28.2477&134.541&123.743\\
set-30-4-2-58&-1044.52&-1044.17&-1044.04&-1042.2&69.5229&168.242&454.528&471.841\\
set-30-4-3-50&-813.949&-748.958&-748.958&-744.291&46.3874&178.628&346.418&285.502\\
\hline  \noalign{\smallskip}
set-30-5-2-60&-828.387&-820.734&-820.734&-818.061&70.6086&177.594&766.148&735.596\\
set-30-5-3-33&-510.902&-494.661&-494.661&-493.585&22.2189&93.7847&247.071&218.586\\
set-30-5-4-46&-520.127&-511.427&-511.346&-509.775&67.1995&283.42&559.463&563.919\\
\hline  \noalign{\smallskip}
set-30-6-1-10&-1027.64&-1023.3&-1023.24&-1021.25&2.27227&11.5146&70.6774&58.1292\\
set-30-6-3-44&-703.572&-702.96&-702.96&-700.314&34.1835&140.288&521.788&530.05\\
set-30-6-4-25&-448.76&-445.707&-445.673&-444.336&14.1767&77.667&185.765&161.619\\
\hline  \noalign{\smallskip}
set-30-7-1-42&-1773.83&-1746.49&-1742.23&-1637.49&35.5206&63.4244&630.688&592.466\\
set-30-7-2-55&-1486.3&-1448.24&-1448.24&-1442.07&63.6699&139.714&983.371&939.227\\
set-30-7-6-43&-194.096&-193.064&-193.064&-191.185&37.2041&265.017&541.029&560.087\\
\hline  \noalign{\smallskip}
set-30-8-1-25&-1659.66&-1531.5&-1531.5&-1515.58&22.1517&25.6327&404.225&311.929\\
set-30-8-2-58&-1010.24&-1009.01&-1008.49&-999.31&70.9503&171.326&1188.2&1258.59\\
set-30-8-6-60&-386.848&-386.538&-386.538&-386.326&76.0659&461.669&1060.4&1013.67\\
\hline  \noalign{\smallskip}
set-30-9-2-60&-969.073&-953.641&-953.641&-949.923&75.0906&179.733&1468.48&1335.38\\
set-30-9-5-30&-273.552&-273.307&-273.307&-272.705&38.512&131.216&421.553&382.539\\
set-30-9-7-58&-282.216&-279.421&-279.421&-279.101&85.9947&579.015&1134.66&1140.53\\
\hline  \noalign{\smallskip}
set-30-10-2-29&-565.335&-563.997&-563.919&-561.784&33.2646&52.5847&557.768&470.508\\
set-30-10-3-31&-506.954&-481.015&-481.015&-478.257&20.9386& 77.4038 &632.39 &542.151\\
set-30-10-8-60&-371.855&-371.216&-371.195&-371.061&87.6211&702.363&1391.21&1329.4\\
\hline
\end{tabular}}
\end{table}
We use the name ``set-$n$-$l$-$k$-$m$" to denote different sets of test problems, where $n$ denotes the dimension of decision variable $x$, $l$ denotes the number of quadratic constraints,
$k$ denotes the number of convex quadratic constraints, and $m$ denotes the number of linear constraints.
We test numerical experiments with $l$ changing from 1 to 10, $k$ changing from 1 to $l-1$ and $m$ changing from 1 to 60, and report  numerical results in Table \ref{tab3} with the examples whose ${\rm(SDP_{GSRT})}$ have  large improvement.

In Table \ref{tab3}, RLT denotes the conic relaxation ${\rm(SDP_{RLT})}$, SOC-RLT denotes the conic relaxation ${\rm(SDP_{ SOC\mbox{-}RLT})}$, GSRT\mbox{-}A denotes the conic relaxation ${\rm(SDP_{GSRT\mbox{-}A})}$ and GSRT-B denotes the conic relaxation ${\rm(SDP_{GSRT\mbox-B})}$, according to their definitions in Section 2.
The number of RLT constraints is $m(m-1)$. The number of SOC-RLT
constraints, which are SOC representable constraints, is $km$. The number of convex quadratic (SOC representable) constraints and that of linear constraints in GSRT constraints are $2(l-k)m+2(l-k)$
and $(l-k)m$, respectively.
Also, to illustrate the effect of the GSRT relaxations, we kick out the examples whose SDP+RLT relaxation are exact, infeasible or unbounded.

We can conclude from Table 3 that  a dominance relationship of $\rm RLT\leq  SOC\mbox{-}RLT\leq GSRT\mbox{-}A\leq GSRT\mbox{-}B$ holds for
the lower bound and a dominance relationship of $\rm RLT\leq SOC\mbox{-}RLT\leq GSRT\mbox{-}B$ or $\rm GSRT\mbox{-}A$ holds for the CPU time.
The tighter lower bounds of both $\rm(SDP_{GSRT\mbox{-}A})$ and $\rm(SDP_{GSRT\mbox{-}B})$   than $\rm(SDP_{ SOC\mbox{-}RLT})$, albeit the increased CPU time cost, are reasonable because of the additional valid inequalities. The comparison
of the lower bounds further shows  an interesting result that the lower bounds of GSRT-B are always better than or equal to the lower bounds of GSRT\mbox{-}A, whose proof remains as an open problem.
For most problem sets,
the CPU time satisfies the following inequality $\rm GSRT\mbox{-}B\leq GSRT\mbox{-}A$.
We also conclude from the table that the number of linear and SOC\ constraints significantly affects the CPU time for different relaxations.
An increment of linear constraints largely increases the number of SOC constraints in SOC-RLT, GSRT\mbox{-}A  and  GSRT-B,  thus increasing
the CPU time significantly. For instances with the same number of quadratic constraints and similar number of linear constraints,
more nonconvex quadratic constraints lead to larger CPU time in GSRT\mbox{-}A and GSRT-B, because a nonconvex quadratic constraint generates
SOC constraints about two times more than a  convex quadratic constraint does and has one more dimension in the lifted matrix.

Since we do not know the optimal value of the examples in Table \ref{tab3}, we could not measure the improvement of the GSRT constraint precisely. In the  following Figures \ref{figeg5}--\ref{figeg15}, we will show that the improvement can be significant for some class of problems. To measure the effect of the GSRT relaxations, we define the improvement ratio as
$${\rm improv.ratio}=\frac{v({\rm SDP_{GSRT}})-v{\rm(SDP_{RLT})}}{ v{\rm(SDP_{RLT})}}. $$
We set the test problems the same as those in Table \ref{tab3} except that the negative eigenvalues in the quadratic constraints have different number of eigenvalues (which is denoted by $\phi$ in Figures \ref{figeg5}--\ref{figeg15}), and $Q_0=I-\sum_{i}^n Q_i$ to ensure the boundedness of the  relaxations. We also set the dimension of the problem as $n=20$, the number of quadratic constraint as $l=5$ and all the quadratic constraints are nonconvex, i.e., $k=0$,  and the linear constraints $m$ changing from 1 to 40.  For each problem setting, we   compute 10 random examples  and illustrate the mean and maximal improvement in the figures.
 From Figures  \ref{figeg5}--\ref{figeg15}, we conclude that the improvement is significant with average improvement up to 9\%, 5\% and 11\% and maximal improvement up to 30\%, 17\% and 36\% for cases that $\phi=5$,  $\phi=10$ and $\phi=15$, respectively.

\section{Concluding remark}
In this paper, we have presented the GSRT  valid inequalities to tighten the SDP relaxations for nonconvex QCQP problems. While the convex relaxations in the current literature lose their effects when dealing with nonconvex quadratic
constraints, we decompose each nonconvex quadratic constraint to two convex quadratic constraints and  develop GSRT constraints based on the
idea of RLT.
Specifically, our GSRT constraints extend the SOC-RLT constraint, by  linearizing  the product of  any pair of linear constraint and SOC constraint derived from nonconvex quadratic constraints.  Enlightened by the
decomposition-approximation method in \cite{zheng2011convex}, we have further proposed a tighter relaxation with extra RLT, SOC-RLT and GSRT generated
by extra valid linear inequality $\alpha_u\geq u^Tx$.
Extending the idea of the GSRT constraints, we have also derived valid inequalities by linearizing the product of any pair of SOC constraints derived from
all quadratic constraints. We have finally extended the Kronecker product constraint to GSOC constraints and demonstrated its
relationship with the previous relaxations.
Promising performance of our numerical tests make us to believe potential applications of our approaches in branch and bound method algorithms for general QCQP problems.

 While we extend the reach of the RLT-like techniques for almost all different types of constraint-pairs, we also examine the dominance relationships among them in order to remove these dominated valid inequalities from consideration. We now summarize the dominance relationships of different relaxations discussed in this paper in the following
Figure \ref{fig3}.\begin{figure}[!htbp]
\centering
\includegraphics[width=0.8\textwidth]{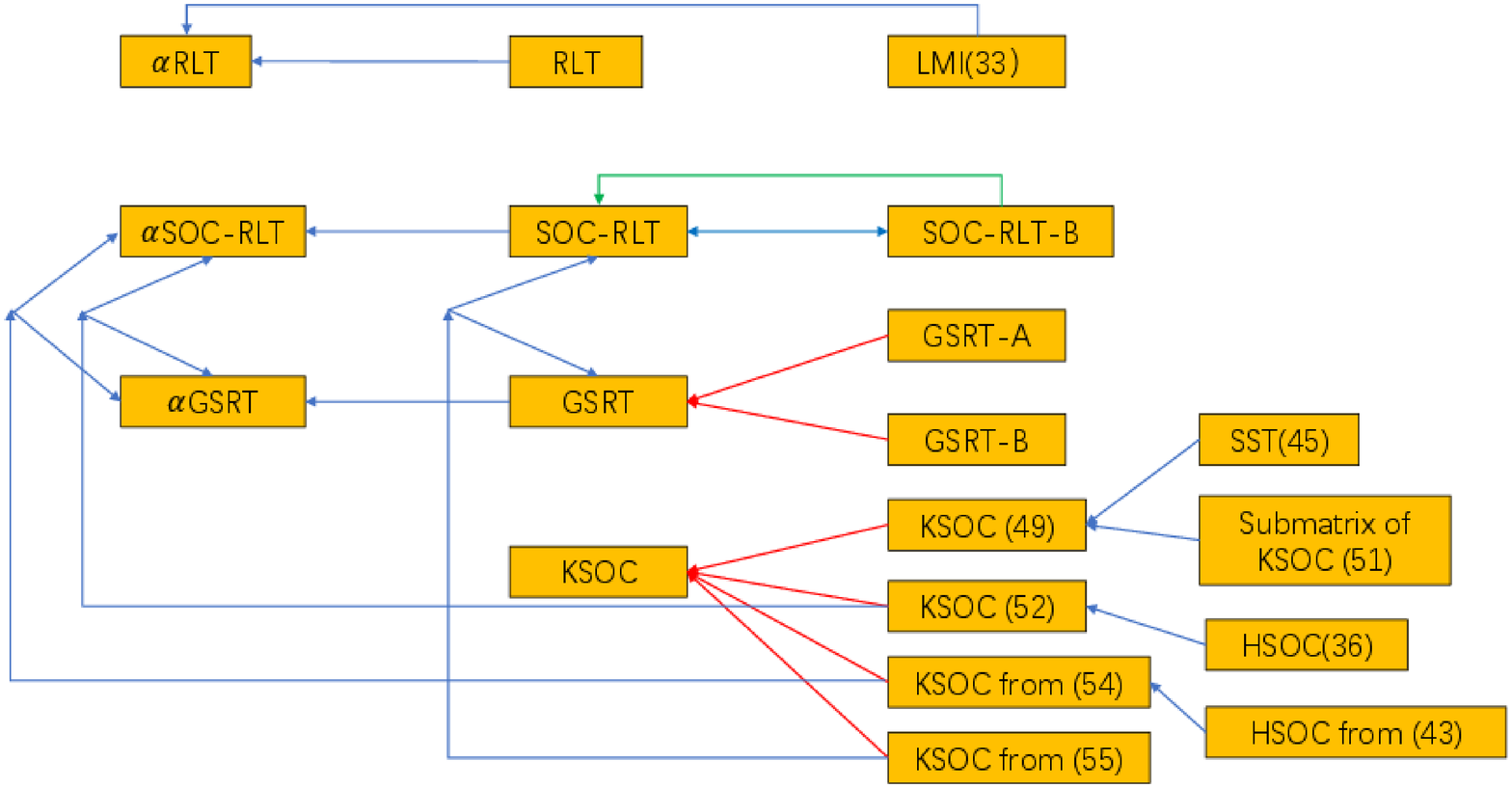}
\caption{This figure shows dominance relationships among different valid inequalities. We use $\rm \alpha RLT$, $\rm \alpha SOC\mbox{-}RLT$ and $\rm \alpha GSRT$ to denote different valid inequalities generated from RLT, SOC-RLT and GSRT with a redundant linear inequality $u^Tx\leq\alpha_u$, respectively. A blue arrow indicates the direction of the dominance. i.e., the valid inequality at the tip of the arrow dominates the valid inequality at the bottom of the arrow, e.g., $\rm\alpha RLT$ dominates $\rm RLT$. A red arrow indicates the direction of an inclusion, i.e., the valid inequality at the tip of the arrow includes the valid inequality at the bottom of the arrow, e.g., GSRT includes GSRT-A and GSRT-B. Also note that KSOC (52) and (54)  are either dominated by $\rm \alpha GSRT$ or $\rm \alpha SOC\mbox{-}RLT$, depending on whether the SOC (that generates (52) and (54)) is derived from convex or nonconvex quadratic constraints.}\label{fig3}
\end{figure}

In fact, we can further rewrite the objective function as $\min \tau$ and add a new constraint $x_0^TQ_0x_0+c_0^Tx\leq\tau$, with a new variable  $\tau$. The original problem is
then equivalent to minimizing $\tau$ and all the techniques developed in this paper can be applied to the new constraint $x_0^TQ_0x_0+c_0Tx\leq\tau$ to achieve a
tighter lower bound.

 An obvious drawback of the relaxations proposed in this paper is their expensive computational cost due to the involved large number of extra SOC
constraints, which is a general challenge in RLT based optimization algorithms, see \cite{anstreicher2009semidefinite,sherali2002enhancing}. One direction to overcome this computational difficulty is to avoid solving SDP problems by using, instead, linear inequalities to approximate the linear matrix constraint $X\succeq xx^{T}$, which
are also called \textit{semidefinite cutting plane} method \cite{sherali2002enhancing} and \cite{qualizza2012linear}. Another important observation is that a large number of  RLT, SOC-RLT and GSRT constraints are inactive at the optimal solution, which inspires us to consider in our future study  the idea of dynamically adding
\textit{semidefinite cutting planes}. More specifically, we can dynamically add some RLT, SOC-RLT and GSRT constraints which are most violated by the current relaxation solution, rather than including all the RLT, SOC-RLT and GSRT constraints in the beginning.
\section*{Acknowledgements}
The authors gratefully acknowledge the support of
Hong Kong Research Grants Council under Grant 14213716.  The second author is also grateful to
the support from Patrick Huen Wing Ming Chair Professorship of Systems Engineering and Engineering Management.

%
%
%





\end{document}